\documentclass[11pt]{amsart}
\usepackage{xcolor}
\usepackage{mathrsfs}
\usepackage{amsmath,amsthm,amssymb,amscd}
\usepackage{booktabs} 
\usepackage{float}
\usepackage{tikz-cd}
\usepackage[hyperfootnotes=false, colorlinks, linkcolor={magenta}, citecolor={magenta}, filecolor={magenta}, urlcolor={magenta}, plainpages=false, pdfpagelabels]{hyperref}
\usepackage{mathtools}
\usepackage[numbers,sort&compress]{natbib} 
\usepackage{etoolbox} 
\apptocmd{\sloppy}{\hbadness 10000\relax}{}{} 
\usepackage[left=2.5cm, right=2.5cm, top=3cm]{geometry}
\allowdisplaybreaks
\usepackage{enumerate}
\usepackage[shortlabels]{enumitem}
\usepackage{bm}
\usepackage{url}
\usepackage[multiple]{footmisc}

\usepackage{colonequals}
\usepackage{longtable} 
\usepackage[ampersand]{easylist}
\usepackage{ltablex,booktabs}

\usepackage{verbatim,delarray,fullpage}
\usepackage{mathtools}
\usepackage{tabularx}
\usepackage{array}
\newcolumntype{C}[1]{>{\centering\let\newline\\\arraybackslash\hspace{0pt}}m{#1}}
\newcolumntype{M}[1]{>{\centering\arraybackslash}m{#1}}
\newcolumntype{N}{@{}m{0pt}@{}}
\usepackage{multirow,caption}
\usepackage{float}
\usepackage{multirow,caption}

\usepackage{stackengine}
\newcommand\xrowht[2][0]{\addstackgap[.5\dimexpr#2\relax]{\vphantom{#1}}}

\theoremstyle{plain}
\newtheorem{Theorem}{Theorem}[section]
\newtheorem{Corollary}[Theorem]{Corollary}
\newtheorem{Lemma}[Theorem]{Lemma}
\newtheorem{Proposition}[Theorem]{Proposition}

\theoremstyle{definition}

\newtheorem{Observation}[Theorem]{Observation}
\newtheorem{Remark}[Theorem]{Remark}

\newtheorem{Question}[Theorem]{Question}
\newtheorem{Example}[Theorem]{Example}

\newcommand{\thmref}[1]{Theorem~\ref{#1}}
\newcommand{\propref}[1]{Proposition~\ref{#1}}
\newcommand{\remref}[1]{Remark~\ref{#1}}
\newcommand{\lemref}[1]{Lemma~\ref{#1}}
\newcommand{\exampref}[1]{Example~\ref{#1}}
\newcommand{\questref}[1]{Question~\ref{#1}}

\newcommand{\sectionref}[1]{Section~\ref{#1}}

\newcommand{\corref}[1]{Corollary~\ref{#1}}

\newcommand{\tabref}[1]{\ref{#1}}

\newcommand{\hs}{\hspace{0.25cm}}
\newcommand{\neghs}{\mkern-12mu}

\DeclareMathOperator{\GL}{GL}
\DeclareMathOperator{\SL}{SL}
\DeclareMathOperator{\Gal}{Gal}

\DeclareMathOperator{\Frob}{Frob}

\DeclareMathOperator{\Aut}{Aut}

\newcommand{\N}{\mathbb{N}}
\newcommand{\Z}{\mathbb{Z}}
\newcommand{\Q}{\mathbb{Q}}

\newcommand{\C}{\mathbb{C}}

\newcommand{\F}{\mathbb{F}}

\newcommand{\nseq}{\mathfrak{n}}

\DeclareMathOperator{\lcm}{lcm}
\DeclareMathOperator{\rad}{rad}

\DeclareMathOperator{\li}{Li}

\newcommand{\Cyc}{\mathcal{C}}
\newcommand{\Div}{\mathcal{D}}

\makeatletter
\renewcommand*{\p@equation}{}
\makeatother

\newcommand{\Eab}{E^{a,b}}
\newcommand{\Eabp}{\widetilde{E^{a,b}_p}}
\DeclareMathOperator{\ST}{ST}

\title{On the average congruence class bias for cyclicity and divisibility of the groups of $\F_p$-points of elliptic curves}
\author[Lee]{Sung Min Lee}
\address{Sung Min Lee, College of Arts and Science, Dakota State University, Madison, SD 57042}
\email{john.lee1@dsu.edu}
\date{\today}

\begin{document}

\begin{abstract}
    In 2009, W. D. Banks and I. E. Shparlinski studied the average densities of primes $p \leq x$ for which the reductions of elliptic curves of small height modulo $p$ satisfy certain algebraic properties, namely cyclicity and divisibility of the number of points by a fixed integer $m$. In this paper, we refine their results by restricting the primes $p$ under consideration to lie in an arithmetic progression $k \bmod{n}$. Furthermore, for a fixed modulus $n$, we investigate statistical biases among the different congruence classes $k \bmod{n}$ of primes satisfying the aforementioned properties.
\end{abstract}

\maketitle

\section{Introduction and statement of results} \label{Introduction}
\subsection{Introduction}\label{s1.1}
Let $E$ be an elliptic curve over $\Q$, and denote its conductor by $N_E$. Let $p$ be a prime of good reduction for $E$, i.e., $p \nmid N_E$. Then the reduction $\widetilde{E}_p$ of $E$ modulo $p$ is an elliptic curve over $\F_p$. For a fixed $E/\Q$, one can study how often $\widetilde{E}_p$ satisfies a given property as $p$ varies. 

For a given property $\mathcal{X}$, we define the (individual) density of primes $p$ for which $\widetilde{E}_p$ satisfies $\mathcal{X}$ as
$$C_E^{\mathcal{X}} \coloneqq \lim_{x\to \infty} \frac{\#\{p \leq x : p \nmid N_E, \widetilde{E}_p \text{ satisfies } \mathcal{X}\} }{\pi(x)},$$
where $\pi(x)$ denotes the prime-counting function, which is asymptotically equivalent to $\li(x) \coloneqq \int_2^x dt/\log t$.

Given a prime $p$ of good reduction for $E$, the $\F_p$-rational points on $\widetilde{E}_p$ form a finite abelian group with at most two generators. That is, there exists positive integers $d_p(E) \mid e_p(E)$ such that
$$\widetilde{E}_p(\F_p) \cong \Z/d_p(E) \Z \oplus \Z/e_p(E) \Z.$$
A natural question is how $\#\widetilde{E}_p(\F_p)$ is distributed as $p$ varies. For instance, for a fixed positive integer $m$, how many primes $p \leq x$ satisfy $\#\widetilde{E}_p(\F_p) = m$? Let $a_p(E)$ denote the trace of Frobenius, which is known to be equal to $p+1-\#\widetilde{E}_p(\F_p)$. By the Hasse-Weil bound, we have
$$|a_p(E)| = |p+1-\#\widetilde{E}_p(\F_p)| \leq 2\sqrt{p}.$$
Thus, there are at most finitely many primes $p$ for which $\#\widetilde{E}_p(\F_p) = m$. In particular, the density is zero.

A slight modification of this question is: given a positive integer $m$, how often is $\#\widetilde{E}_p(\F_p)$ divisible by $m$? We call such primes $p$ \textbf{primes of $m$-divisibility} for $E$. In other words, we seek to understand the asymptotic behavior of 
$$\pi_E^{\Div(m)}(x) \coloneqq \#\{p \leq x : p \nmid mN_E, m \mid \#\widetilde{E}_p(\F_p)\}, \hs \text{ as } x \to \infty.$$
A.\ C.\ Cojocaru \cite{MR2076566} proved that there exists a constant $C^{\Div(m)}_E$ for which
$$\pi_E^{\Div(m)}(x) \sim C^{\Div(m)}_E \cdot \frac{x}{\log x}, \hs \text{ as } x\to \infty.$$
The constant $C^{\Div(m)}_E$ can be computed using Galois representations and the Chebotarev density theorem as follows.

Let $G_\Q$ be the absolute Galois group, and let $E[m] \subseteq E(\overline{\Q})$ be the subgroup of $m$-torsion points of $E$. Define $\Q(E[m])$, the $m$-division field, as the field obtained by adjoining the coordinates of $P \in E[m]$. Since $G_\Q$ acts on $E[m]$, upon choosing a $\Z/m\Z$-basis, we define the mod $m$ Galois representation:
$$\rho_{E,m} : G_\Q \to \Aut(E[m]) \cong \GL_2(\Z/m\Z).$$
For a prime $p \nmid mN_E$, let $\text{Frob}_p$ be a Frobenius element in $G_\Q$ lying over $p$. Then we have
\begin{equation}\label{easyfact}
    \det(\rho_{E,m}(\text{Frob}_p)) \equiv p \mkern-12mu \pmod m \hs \text{ and } \hs \text{Tr}(\rho_{E,m}(\text{Frob}_p)) \equiv a_p(E) \mkern-12mu  \pmod m.
\end{equation}
Define
\begin{equation}\label{CmEdefinition}
    D_m(E) \coloneqq \{ g \in \Gal(\Q(E[m])/\Q) : \det (g) + 1 - \text{Tr}(g) \equiv 0 \mkern-12mu \pmod m\}.
\end{equation}
From \eqref{easyfact} and \eqref{CmEdefinition}, it follows that
$$m \mid \#\widetilde{E}_p(\F_p) \iff  \rho_{E,m}(\text{Frob}_p) \in D_m(E).$$
Since $D_m(E)$ forms a conjugacy class in $\Gal(\Q(E[m])/\Q)$, the right-hand side holds for any choice of $\Frob_p$ lying above $p$. By the Chebotarev density theorem, we obtain
$$C_E^{\Div(m)} = \frac{\# D_m(E)}{[\Q(E[m]):\Q]}.$$
Additionally, Cojocaru showed that the density is approximately $1/m$.

Another natural question is how $d_p(E)$ is distributed as $p$ varies. In particular, how often do we have $d_p(E) = 1$, i.e., when $\widetilde{E}_p(\F_p)$ is cyclic? We refer to such primes $p$ as \textbf{primes of cyclic reduction} for $E$. Reformulating the question, we seek to understand the asymptotic behavior of
$$\pi^{\Cyc}_E(x) \coloneqq \#\{p \leq x : p \nmid N_E, d_p(E) = 1\}, \hs \text{ as } x\to \infty.$$


The first major progress on this problem was made by J.-P. Serre, who observed that it shares similarities with Artin's primitive root conjecture. This conjecture was conditionally proven under Generalized Riemann Hypothesis (GRH) by C.\ Hooley \cite{MR207630}. Serre's result states that under GRH, there exists a constant $C^{\Cyc}_E$ for which
$$\pi_E^{\Cyc}(x) \sim C_E^{\Cyc} \cdot \frac{x}{\log x}, \hs \text{ as } x \to \infty,$$
where
$$C_E^{\Cyc} = \sum_{n \geq 1} \frac{\mu(n)}{[\Q(E[n]):\Q]}.$$
(If the constant vanishes, we interpret this meaning that there are only finitely many primes of cyclic reduction for $E$.) Moreover, he observed that $C_E^{\Cyc} > 0$ if and only if $\Q(E[2]) \neq \Q$. (The proof can be found in \cite[pp.160-161]{MR0698163}.)

Following Serre's result, significant progress has been made. In 1983, M. Ram Murty removed GRH for CM elliptic curves \cite{MR0698163}. In 1990, R. Gupta and M. Ram Murty \cite{MR1055716} established that if $\Q(E[2]) \neq \Q$, then $\pi_E^{\Cyc}(x) \gg_E x/\log^2 x$. In 2002, Cojocaru established Serre's results under a weaker hypothesis known as quasi-GRH and also improved the error bound \cite{MR1932460}. The error term has been extensively analyzed under various assumptions; details can be found in \cite{MR4033729}. Recently, F. Campagna and P. Stevenhagen extended the GRH-conditional results to arbitrary number fields $K$ and elliptic curves $E/K$ \cite{campagna2022cyclic}.

The density $C^{\Cyc}_E$ also has been studied in various papers. In 2009, N.\ Jones explicitly calculated $C^{\Cyc}_E$ for Serre curves $E$. Later, J.\ Brau refined this result for elliptic curves with abelian entanglements \cite{brau2017character}. Also, J.\ Mayle and Rakvi calculated the density for a family of elliptic curves known as relative Serre curves \cite{mayle2023serre}.

In 2009, W.\ D.\ Banks and I.\ E.\ Shparlinski studied the average version of these problems. Any elliptic curve over $\Q$ can be expressed as
$$ \Eab : Y^2 = X^3+aX+b,$$
for some $a,b \in \Z$ and $4a^3+27b^2 \neq 0$. We define the height of $E^{a,b}$ as $\max\{|a|^3,|b|^2\}$. Define
\begin{align*}
    \pi^{\Cyc}_{a,b}(x) &\coloneqq \#\{p \leq x : p \nmid N_{\Eab},d_p(\Eab) = 1\}, \\
    \pi^{\Div(m)}_{a,b}(x) &\coloneqq \#\{p \leq x : p \nmid mN_{\Eab}, m \mid \#\Eabp(\F_p)\}.
\end{align*}
If $4a^3+27b^2 = 0$, we define both $\pi^{\Cyc}_{a,b}(x)$ and $\pi^{\Div(m)}_{a,b}(x)$ to be the zero function. Two of the main theorems from \cite{MR2570668} are summarized as follows.
\begin{Theorem}[\protect{\cite[Theorems 18 and 19]{MR2570668}}]\label{BSsummary}
Let $\epsilon > 0$ and $K > 0$ be fixed. Let $A \coloneqq A(x)$ and $B \coloneqq B(x)$ be integers satisfying
\begin{equation}\label{conditionsonAB}
    x^\epsilon \leq A, B \leq x^{1-\epsilon}, \hs x^{1+\epsilon} \leq AB.
\end{equation}
\begin{itemize}
    \item There exists a positive constant $C^{\Cyc}$ for which
    $$\frac{1}{4AB}\sum_{|a| \leq A}\sum_{|b| \leq B} \pi^{\Cyc}_{a,b}(x) = \left(C^{\Cyc} + O_{\epsilon, K}\left( \log^{-K}x\right)\right)\li(x), \hs \text{ as } x \to \infty.$$
    \item There exists a positive constant $C^{\Div(m)}$ such that for $m \leq \log^K x$,
    $$\frac{1}{4AB} \sum_{|a| \leq A} \sum_{|b| \leq B} \pi^{\Div(m)}_{a,b}(x) = \left( C^{\Div(m)}  + O_{\epsilon, K} \left(\log^{-K} x \right)\right)\li(x), \hs \text{ as } x \to \infty.$$
\end{itemize}
\end{Theorem}
Note that the condition \eqref{conditionsonAB} becomes empty if $\epsilon > 1/3$. Hence, throughout this paper, we assume $\epsilon < 1/3$.

\begin{Remark}\label{rem}
Let $A$ and $B$ satisfy the conditions listed in \thmref{BSsummary}. Consider the quantity
$$\frac{1}{4AB} \sum_{|a| \leq A} \sum_{|b| \leq B} \pi_{a,b}^{\mathcal{X}}(x).$$
This expression counts the average number of primes $p \leq x$ for which $\widetilde{E}_p$ satisfies $\mathcal{X}$ for approximately $4AB$ different models of elliptic curves $E/\Q$. We define the \textbf{average density} of such primes as 
$$C^{\mathcal{X}} \coloneqq \lim_{x\to \infty} \frac{1}{4AB \pi(x)} \sum_{|a| \leq A} \sum_{|b| \leq B} \pi^{\mathcal{X}}_{a,b}(x).$$
Thus, we refer to $C^{\Cyc}$ and $C^{\Div(m)}$ as the average densities of primes of cyclic reduction and $m$-divisibility, respectively.
\end{Remark}

In this paper, we extend \thmref{BSsummary} by incorporating a congruence class condition and investigating whether a statistical bias exists in the distribution of these primes across different congruence classes. Fix a positive integer $n$ and let $k$ be coprime to $n$. Define
\begin{align*}
    \pi^{\Cyc}_{a,b}(x; n,k) &\coloneqq \#\{p \leq x : p \nmid N_{\Eab}, \, d_p(\Eab) = 1, \, p \equiv k \mkern-12mu \pmod n\}, \\
    \pi^{\Div(m)}_{a,b}(x; n,k) &\coloneqq \#\{p \leq x : p \nmid mN_{\Eab}, \, m \mid \#\Eabp(\F_p), \, p \equiv k \mkern-12mu \pmod n\}.
\end{align*}
As before, if $4a^3+27b^2$ vanishes, we regard both of them as the zero function. We now state the main theorems of this paper.
\begin{Theorem}\label{maintheorem}
    Let $x > 0$, $\epsilon > 0$, and $K > 0$ be fixed. Let $A \coloneqq A(x)$ and $B \coloneqq B(x)$ satisfy
$$x^{\epsilon} \leq A , B \leq x^{1-\epsilon}, \hs  
 \hs x^{1+\epsilon} \leq AB.$$
 \begin{itemize}
     \item Fix $\sigma < \epsilon/(16+4\epsilon)$. For any $n \leq \log^K x$ and $k$ coprime to $n$, there exists a positive constant $C^{\Cyc}_{n,k}$ such that
    $$\frac{1}{4AB} \sum_{|a| \leq A} \sum_{|b| \leq B} \pi_{a,b}^{\Cyc}(x;n,k) = \left(C^{\Cyc}_{n,k} + O_{\epsilon, \sigma, K}\left( \frac{\log x}{nx^{\sigma}} \right)\right) \li(x), \hs \text{ as } x\to \infty.$$

    \item  Choose positive integers $n$ and $m$ such that
    $$\lcm\left(n,  \prod_{q^j \parallel m} q^{\lceil j/2\rceil} \right) \leq \log^{K+1}x.$$
    For any $k$ coprime to $n$, there exists a positive constant $C^{\Div(m)}_{n,k}$ such that
    $$\frac{1}{4AB} \sum_{|a| \leq A} \sum_{|b| \leq B} \pi_{a,b}^{\Div(m)}(x;n,k) = \left(C^{\Div(m)}_{n,k} + O_{\epsilon, K}\left(\frac{1}{n \log^K x} \right) \right)\li(x), \hs \text{ as } x\to \infty.$$
 \end{itemize} 
\end{Theorem}
Explicit formulas for $C^{\Cyc}_{n,k}$ and $C^{\Div(m)}_{n,k}$ are given in \sectionref{s3.1}.

\begin{Remark}\label{Sato-Tate}
    Let $E/\Q$ be a non-CM elliptic curve of conductor $N_E$ and $p \nmid N_E$ be a prime. The Hasse-Weil bound states that $|a_p(E)| \leq 2\sqrt{p}$. Thus, there exists a real number $\psi_E(p)$, called the Frobenius angle of $E$ at $p$, in $[0,\pi]$ for which
    \begin{equation}\label{Frobangle}
        a_p(E) = 2\sqrt{p} \cos (\psi_E(p)).
    \end{equation}
    The famous Sato-Tate conjecture asserts that, given $0 \leq \alpha < \beta \leq \pi$,
    \begin{equation}\label{SatoTatemeasure}
              \# \{p \leq x : p \nmid N_E, \psi_E(p) \in [\alpha,\beta]\} \sim \frac{2}{\pi} \int_\alpha^\beta \sin^2 \theta d\theta =: \mu_{\ST}(\alpha,\beta).
    \end{equation}
    Recently, L.\ Clozel, M.\ Harris, N.\ Shpherd-Barron, and R.\ Taylor proved that the conjecture holds for elliptic curves over totally real field with multiplicative reduction at some prime. (See \cite{MR2470688} and  references therein.)

    Banks and Shparlinski proved the Sato-Tate conjecture on average \cite[Theorem 17]{MR2570668}. Set $\epsilon > 0$, and let $A$ and $B$ satisfy \eqref{conditionsonAB}. Then, there exists $\delta > 0$ for which
    $$\frac{1}{4AB} \sum_{|a| \leq A} \sum_{|b| \leq B} \#\{ p \leq x : p \nmid N_{\Eab}, \psi_{\Eab}(p) \in [\alpha,\beta] \} = \left(\mu_{\ST}(\alpha,\beta) + O_{\epsilon}(x^{-\delta})\right) \pi(x), \hs \text{ as }x\to \infty.$$
    In \remref{STACC}, we briefly sketch out the proof of the average Sato-Tate conjecture under the arithmetic progression condition.
\end{Remark}
By a similar reasoning as in \remref{rem}, the quantities in \thmref{maintheorem} can be understood as the average densities of primes of cyclic reduction and $m$-divisibility within an arithmetic progression.

We say the primes of property $\mathcal{X}$ exhibit an \textbf{average congruence class bias} modulo $n$ if there exists $k$ and $k'$ coprime to $n$ for which $C^{\mathcal{X}}_{n,k} \neq C^{\mathcal{X}}_{n,k'}$. The following proposition clarifies whether primes of cyclic reduction and $m$-divisibility display such biases.
\begin{Proposition}\label{Densities}
Fix a positive integer $n \geq 2$ and let $k$ be coprime to $n$. Then we have
\begin{itemize}
    \item $C^{\Cyc}_{n,1} \leq C^{\Cyc}_{n,k} \leq C^{\Cyc}_{n,-1}$. If $n$ is a power of two, then $C^{\Cyc}_{n,1} = C^{\Cyc}_{n,-1}$. Otherwise, $C^{\Cyc}_{n,1} < C^{\Cyc}_{n,-1}$. In particular, primes of cyclic reduction always have an average congruence class bias modulo $n$, except when $n$ is a power of two.
    \item $C^{\Div(m)}_{n,k} = C^{\Div(m)}/\phi(n)$ if one of the following holds:
    \begin{itemize}
    \item $n$ and $m$ are coprime.
    \item $m \in \{2,4\}$.
    \end{itemize}
    In particular, primes of $m$-divisibility do not have an average congruence class bias modulo $n$ either if $n$ and $m$ are coprime or if $m\in \{2,4\}$.
\end{itemize}
\end{Proposition}
\begin{Remark}   
    \propref{Densities} reveals that $C^{\Cyc}_{n,k}$ follows a universal pattern independent of $n$. It also states that when $n$ is a power of two, there is no average congruence class bias modulo $n$. However, this does not hold at the individual level. In particular, Jones and the author \cite[Example 1.4]{jones2023acyclicity} provide an example of an elliptic curve $E/\Q$ with infinitely many primes of cyclic reduction, none of which are congruent to 3 or 5 modulo 8. (Details can be found in \sectionref{s2.2}.)

    The second part of the proposition asserts that primes of $m$-divisibility do not exhibit an average congruence class bias modulo $n$ when $n$ and $m$ are coprime or $m \in \{2,4\}$. However, this condition is not necessary; for instance, primes of $8$-divisibility exhibit no average congruence class bias modulo 6. (See Table \tabref{Table3}.) 

    This raises a natural question: Are the primes of $m$-divisibility for an individual elliptic curve $E$ evenly distributed modulo $n$ when $n$ and $m$ are coprime? The answer is no. In \sectionref{s2.1}, we show that the coprimality of $n$ and $m$ alone does not guarantee an equal distribution of primes of $m$-divisibility for $E$ modulo $n$.
\end{Remark}
Let us conclude this section with a number of possible future research directions. In 1999, C.\ David and F.\ Pappalardi proved the Lang-Trotter conjecture on average \cite{MR1677267}. (Their error term was later refined by S.\ Baier \cite{MR2376806}.) Thus, a natural follow up is to establish the asymptotic formula for
$$\frac{1}{4AB} \sum_{|a| \leq A} \sum_{|b| \leq B} \#\{p \leq x : p \nmid N_{\Eab}, \, a_p(\Eab) = t,\, p \equiv k \mkern-12mu \pmod n\},$$
and determine whether there is an average congruence class bias for primes of a fixed Frobenius trace. This study would complement the work of S.\ Jarov, A.\ Khandra, and N.\ Walji who examined congruence class bias at different moduli $n$ for a specific family of elliptic curves \cite{JKW}.

Inspired by $m$-divisibility, one can investigate the average densities of primes for which $\widetilde{E}_p(\F_p)$ or $d_p(E)$ satisfies certain conditions. For instance, as $p$ varies,
\begin{itemize}
    \item how often is $d_p(E) = m$?
    \item how often is $\widetilde{E}_p(\F_p) = \langle \bar{a} \pmod p \rangle$, where $a \in E(\Q)$ is a point of infinite order?
    \item how often is $\#\widetilde{E}_p(\F_p)$ square-free?
    \item how often is $\#\widetilde{E}_p(\F_p)$ prime?
\end{itemize}
Cojocaru compiled results on the individual densities of these primes in \cite{MR2076566}. One could employ the method of Banks and Shparlinski to study their average densities. Furthermore, adding a congruence class condition could reveal whether these primes exhibit average congruence class biases.

\subsection{Contents of the Paper and Notation} 
This paper is structured as follows. In \sectionref{s2}, we examine the congruence class bias of primes of $m$-divisibility and cyclic reduction for an individual elliptic curve. We also compare the individual and the average results. In \sectionref{s3}, we express the average densities $C^{\mathcal{X}}$ and $C^{\mathcal{X}}_{n,k}$ explicitly and adapt the method of Banks and Shparlinski to incorporate an arithmetic progression. In \sectionref{s4}, we define the mean behavior of an arithmetic function. Modifying techniques from \cite[Sections 2 and 3]{MR2116969}, we prove key lemmas for the main results. In \sectionref{s5}, we present the proofs of \propref{Densities} and \thmref{maintheorem}. Finally, \sectionref{s6} provides tables of values of $C^{\mathcal{X}}_{n,k}$ for small $n$.

Below is a list of the notations used in this paper.
\begin{itemize}
    \item Given two functions $f$ and $g$ of $x$, we write $f \ll g$ or $f = O(g)$ if there exists a constant $C > 0$ and $x_0$ such that $|f(x)| \leq C|g(x)|$ for any $x > x_0$. If the constant $C$ depends on a parameter $m$, we denote $f \ll_m g$ or $f = O_m(g)$.
     \item We denote $f \sim g$ if $\lim_{x\to \infty} \frac{f(x)}{g(x)} = 1$. Also, $f = o(g)$ if $\lim_{x\to \infty} \frac{f(x)}{g(x)} = 0$.
    \item $p, q, \ell, \ell'$ denote primes. $n$ denotes a positive integer (usually greater than 1) and $k$ an integer coprime to $n$ or a congruence class in $(\Z/n\Z)^\times$, depending on context.
    \item $\prod_q$ denotes that a product taken over all primes $q$. $\sum_{p \equiv k (n)}$ denotes a sum taken over all primes $p$ such that $p \equiv k \pmod n$. 
    \item $\li(x)$ denotes the logarithmic integral function, i.e.,
    $$\li(x) \coloneqq \int_2^x \frac{dt}{\log t}.$$
    It is well-known that $\li(x) \sim x \log^{-1}x \sim \pi(x)$.
\end{itemize}

\subsection{Acknowledgement} I am grateful to my advisor, Nathan Jones, for his insightful guidance and helpful comments. I also thank him for providing \exampref{5div} and \exampref{nontriv}. I would like to thank Jacob Mayle and Tian Wang for carefully reading the manuscript and providing detailed feedback. I appreciate Igor E. Shparlinski for his advice and the anonymous referees for their thorough comments and valuable suggestions for improvement. Finally, I thank the anonymous referees for their valuable comments and suggestions.

\section{Congruence Class Bias for an Individual Elliptic Curve}\label{s2}
\subsection{Primes of \texorpdfstring{$m$}{m}-divisibility in an arithmetic progression}\label{s2.1}
Define
$$C^{\mathcal{X}}_{E,n,k} \coloneqq \lim_{x \to \infty} \frac{\#\{p \leq x : p \nmid N_E, \widetilde{E}_p \text{ satisfies } \mathcal{X}, p \equiv k \mkern-12mu \pmod n\}}{\pi(x)}.$$
We say that the primes of property $\mathcal{X}$ for $E$ exhibit a congruence class bias modulo $n$ if there exists $k$ and $k'$ coprime to $n$ such that $C^{\mathcal{X}}_{E,n,k} \neq C^{\mathcal{X}}_{E,n,k'}$.

Let us first examine whether primes of $m$-divisibility for $E$ show such a bias. Let $\zeta_n$ be a primitive $n$-th root of unity, and let $\sigma_k \in \Gal(\Q(\zeta_n)/\Q)$ be the automorphism sending $\zeta_n \mapsto \zeta_n^k$. Consider restriction maps
\begin{align*}
    \rho_1 : & \Gal(\Q(E[m])/\Q) \to \Gal(\Q(E[m])\cap \Q(\zeta_n)/\Q), \\
    \rho_2 : & \Gal(\Q(\zeta_n)/\Q) \to \Gal(\Q(E[m])\cap \Q(\zeta_n)/\Q).
\end{align*}
Then, $\Gal(\Q(E[m])\Q(\zeta_n)/\Q)$ can be realized as a fiber product of $\Gal(\Q(E[m])/\Q)$ and $\Gal(\Q(\zeta_n)/\Q)$ along $\rho \coloneqq (\rho_1,\rho_2)$, yielding the following diagram:
\begin{center}
    \begin{tikzcd}
\Gal(\Q(E[m])/\Q) \times_\rho \Gal(\Q(\zeta_n)/\Q) \arrow[r, "\pi_2"] \arrow[d, "\pi_1"] & \Gal(\Q(\zeta_n)/\Q) \arrow[d, "\rho_2"] \\
    \Gal(\Q(E[m])/\Q) \arrow[r, "\rho_1"] & \Gal(\Q(E[m]) \cap \Q(\zeta_n)/\Q).
    \end{tikzcd}
\end{center}
(Here, $\pi_1$ and $\pi_2$ denote the natural projections.)

For $p \nmid nmN_E$, let $\text{Frob}_p$ be a Frobenius element in $\Gal(\Q(E[m])\Q(\zeta_n)/\Q)$ above $p$. Let $D_m(E)$ be as in \eqref{CmEdefinition}. Define
$$D_{m}^{n,k}(E) \coloneqq\{g \in \Gal(\Q(E[m])\Q(\zeta_n)/\Q) : \pi_1(g) \in D_m(E), \pi_2(g) = \sigma_k\}.$$
Then we have
$$p \equiv k \mkern-12mu  \pmod n \text{ and } m \mid \#\widetilde{E}_p(\F_p) \iff \text{Frob}_p \in D_{m}^{n,k}(E).$$
Since $D_m^{n,k}(E)$ forms a conjugacy class, this holds for any choice of $\Frob_p$ above $p$. Thus, by the Chebotarev density theorem, we have
$$C_{E,n,k}^{\Div(m)} = \frac{\#D_{m}^{n,k}(E)}{[\Q(E[m])\Q(\zeta_n) : \Q]}.$$
If $\Q(E[m])$ and $\Q(\zeta_n)$ intersect trivially, then $\Gal(\Q(E[m])\Q(\zeta_n)/\Q)$ is simply the direct product of $\Gal(\Q(E[m])/\Q)$ and $\Gal(\Q(\zeta_n)/\Q)$. In this case,
$$C^{\Div(m)}_{E,n,k} = \frac{\#D_{m}^{n,k}(E)}{[\Q(E[m])\Q(\zeta_n):\Q]} = \frac{\#D_m(E)}{\phi(n) [\Q(E[m]):\Q]}.$$
Thus, primes of $m$-divisibility for $E$ exhibit no congruence class bias modulo $n$ if $\Q(E[m])\cap \Q(\zeta_n) = \Q$.

\propref{Densities} states that there is no average congruence class bias modulo $n$ for primes of $m$-divisibility, provided $(m,n) = 1$. This is a weaker condition than $\Q(E[m]) \cap \Q(\zeta_n) = \Q$, since $\Q(\zeta_m) \subset \Q(E[m])$ via the Weil pairing. One might then ask: what happens if $\Q(E[m]) \cap \Q(\zeta_n) \neq \Q$?

This problem, known as the missing Frobenius trace problem, was first explored by N.\ Jones and K.\ Visseut \cite{jones2021elliptic}, leveraging the comprehensive database compiled by D.\ Zywina \cite{zywina2015possible}. Let $E$ be an elliptic curve with the property that there exists $a \in \Z \setminus \{1\}$ such that for any prime $p$ of good reduction, $a_p(E) \neq a$. Choose $m$ coprime to $a-1$. Since $\#\widetilde{E}_p(\F_p) = p+1-a_p(E)$, we have
$$m \mid \#\widetilde{E}_p(\F_p) \iff a_p(E) \equiv p+1 \pmod m.$$
Taking $n = m$ and $p \equiv a-1 \pmod n$, we see that $p$ cannot be a prime of $m$-divisibility for $E$ since $a_p(E)$ always avoids $a$. If there exists another congruence class $k \pmod n$ with $D_m^{n,k}(E) \neq \emptyset$, then
$$C^{\Div(m)}_{E,n,a-1} = 0 \hs \text{ while } \hs C^{\Div(m)}_{E,n,k} > 0.$$
Therefore, the primes of $m$-divisibility for such $E$ exhibit a congruence class bias modulo $n$.

Jones and Vissuet \cite[Section 5.4]{jones2021elliptic} present $j$-invariants and models of elliptic curves whose mod 5 Galois images are conjugate to
$$G_5 = \left\{ \begin{pmatrix}
    a^2 & * \\ 0 & a
\end{pmatrix}: a \in (\Z/5\Z)^\times\right\}.$$
For such curves $E$, $a_p(E)$ must be congruent to 0, 1, or 2 modulo 5 for all good prime. Thus, $D_5^{5,k}(E)$ is nonempty if $k \in \{1,4\}$, and empty otherwise. Below is an example, kindly provided by Jones, obtained by plugging $t = 1$ in their model $\mathcal{E}_{5,1,2}$.
\begin{Example}\label{5div}
Consider an elliptic curve given by
$$E: Y^2 + Y = X^3+X^2-8X + 19.$$
This curve has conductor 275, and its LMFDB label is 275.b3 \cite{lmfdb}. The primes of $5$-divisibility for $E$ exhibit a congruence class bias modulo $5$. In particular, $C^{\Div(5)}_{E,5,k}$ is zero if $k \in \{2,3\}$ and positive if $k \in \{1,4\}$.
\end{Example}

Here is another example, also provided by Jones, with $n$ and $m$ being coprime.
\begin{Example}\label{nontriv}
    Consider the elliptic curve given by
    $$E : Y^2 = X^3+X^2-2X-1.$$
    This curve has conductor 784, and its LMFDB label is 784.g2. We have
    $$C^{\Div(2)}_{E,7,k} = \begin{cases} 1/6 &\text{ if } k \in \{1,6\} \\0 &\text{ otherwise} \end{cases}$$
\end{Example}

\begin{proof}
    One can verify that $\Q(E[2])$ is a totally real subfield of $\Q(\zeta_7)$ with $[\Q(E[2]):\Q] = 3$. Thus, $G_2(E) \coloneqq \rho_{E,2}(G_\Q)$ is the unique subgroup of $\GL_2(\Z/2\Z)$ of order 3, i.e.,
    $$\Gal(\Q(E[2])/\Q) \cong \left\{ \begin{pmatrix} 1 & 0 \\ 0 & 1\end{pmatrix}, \begin{pmatrix} 1 & 1 \\ 1 & 0\end{pmatrix}, \begin{pmatrix} 0 & 1 \\ 1 & 1\end{pmatrix} \right\},$$
    where only the identity element has an even trace.
    
    Given an odd prime $p \nmid N_E$, consider the Frobenius element in $\Gal(\Q(\zeta_7)/\Q)$ at $p$. Let us denote its restriction to $\Q(E[2])$ by $\text{Frob}_p$. If $p \equiv \pm 1 \pmod 7$, then $\text{Frob}_p$ is the identity element and has an even trace. Otherwise, $\Frob_p$ has an odd trace. Since
    $$\#\widetilde{E}_p(\F_p) = p+1-\text{Tr}(\Frob_p) \equiv \begin{cases}
        1 \pmod 2 & \text{ if } \text{Frob}_p \text{ has an odd trace,} \\
        0 \pmod 2 & \text{ if }\text{Frob}_p \text{ has an even trace,} 
    \end{cases}$$
    we conclude that
    \begin{align*}
        p \equiv \pm 1 \pmod 7 \hs &\Longrightarrow \hs 2 \mid \#\widetilde{E}_p(\F_p), \\
        p \not \equiv \pm 1 \pmod 7 \hs &\Longrightarrow \hs 2 \nmid \#\widetilde{E}_p(\F_p).
    \end{align*}
    This completes the proof.
\end{proof}

We observe that primes of $m$-divisibility for an individual elliptic curve can have a congruence class bias modulo $n$ when $\Q(E[m])\cap \Q(\zeta_n) \neq \Q$. A natural question is how common this phenomenon is. We show that for any Serre curve $E/\Q$, there exists such an $n$ with $m = 2$. Since almost all elliptic curves are known to be Serre curves when ordered by height \cite{MR2563740}, we can conclude that such examples are extremely common.
\begin{Observation}
Let $E/\Q$ be a Serre curve defined by
$$E : Y^2 = X^3 + aX + b.$$
Define its discriminant as $\Delta_E \coloneqq -16(4a^3+27b^2)$, and let $D$ be the square-free part of $\Delta_E$. Set $d = |D|$ if $D \equiv 1 \pmod 4$ and $d = 4|D|$ otherwise.

Since $E$ is a Serre curve, $E$ has Serre entanglement: $\Q(E[2]) \cap \Q(\zeta_n) \supset \Q(\sqrt{D})$. Moreover, $E$ has a surjective mod 2 Galois representation, and hence $\Gal(\Q(E[2])/\Q) \cong \GL_2(\Z/2\Z) \cong S_3$, the symmetric group of three elements.

Now, consider an odd prime $p \nmid N_E$ with $\left( \frac{D}{p} \right) = -1$. Let $\text{Frob}_p$ be a Frobenius element in $\Gal(\Q(E[2])/\Q)$ at $p$. Since $\Frob_p|_{\Q(\sqrt{D})}$ is nontrivial, $\text{Frob}_p$ must lie in the nontrivial coset $gA_3$, which can be identified as
$$gA_3 = \left\{\begin{pmatrix}
    1 & 1 \\ 0 & 1
\end{pmatrix}, \begin{pmatrix}
    1 & 0 \\ 1 & 1
\end{pmatrix}, \begin{pmatrix}
    0 & 1 \\ 1 & 0
\end{pmatrix} \right\}.$$
Each element in $gA_3$ has an even trace, implying $2 \mid \#\widetilde{E}_p(\F_p)$, by the same argument used in \exampref{nontriv}. 

On the other hand, for an odd prime $q \nmid N_E$ with $\left(\frac{D}{q}\right) = 1$ and $\text{Frob}_q \neq 1$, we find that $\text{Frob}_q$ belongs to the trivial coset $A_3$. This means it has an odd trace, leading to $2 \nmid \widetilde{E}_q(\F_q)$. These two elements of odd trace form a conjugacy class in $\GL_2(\Z/2\Z)$. Thus, the Chebotarev density theorem ensures that one-third of primes have odd Frobenius traces.

Let $K$ and $K'$ denote the sets of quadratic non-residues and residues modulo $n$, respectively, each of cardinality $\phi(n)/2$. Then, we obtain
\begin{align*}
    \sum_{k \in K} C_{E,n,k}^{\Div(2)} = \frac{1}{2} \hs \text{ while } \hs
    \sum_{k' \in K'} C_{E,n,k'}^{\Div(2)} = \frac{1}{6}.
\end{align*}
Therefore, primes of $2$-divisibility for $E$ exhibit a congruence class bias.
\end{Observation}

\subsection{Primes of cyclic reduction in an arithmetic progression} \label{s2.2}
This section presents GRH-conditional results by Y.\ Akbal and A.\ M.\ G\"{u}lo\v{g}lu's on the asymptotic behavior of $\pi_E^{\Cyc}(x;n,k)$ \cite{MR4504664}. Define
$$\gamma_{n,k}(\Q(E[d])) \coloneqq \begin{cases} 1 & \text{ if } \sigma_k \text{ fixes } \Q(E[d])\cap \Q(\zeta_n) \text{ pointwise,} \\ 0 &\text{ otherwise.} \end{cases}$$
Now set
\begin{equation}\label{CCEnk}
C^{\Cyc}_{E,n,k} \coloneqq \sum_{d\geq 1} \frac{\mu(d)\gamma_{n,k}(\Q(E[d]))}{[\Q(E[d])\Q(\zeta_n):\Q]},
\end{equation}
where $\mu$ is the M\"obius function.
\begin{Theorem}[\protect{\cite[Theorems 3 and 5]{MR4504664}}]
Let $E/\Q$ be an elliptic curve. If $E$ has CM, assume it has CM by the full ring of integers of an imaginary quadratic field. Fix a positive integer $n$. Assume that GRH holds for the Dedekind zeta functions of the field $\Q(E[d])\Q(\zeta_n)$ for any square-free integer $d\geq 1$. For any $k \in \Z$ with $(k,n) = 1$, we have
$$\pi^{\Cyc}_E(x;n,k) \sim C^{\Cyc}_{E,n,k}\cdot \frac{x}{\log x}, \hs \text{ as } x\to \infty.$$
\end{Theorem}
P.-J. Wong later confirmed that this results hold unconditionally for CM elliptic curves \cite{MR4765791}. However, formula \eqref{CCEnk} does not indicate for which values of $k$ maximizes or minimizes the density $C^{\Cyc}_{E,n,k}$.

In his doctoral thesis, J. Brau \cite{Brau} derived an almost Euler product formula for $C^{\Cyc}_{E,n,k}$ for Serre curves, except for a thin family. Recently, in collaboration with J. Mayle and T. Wang, the author \cite{lee2024opposingaveragecongruenceclass}, extended Brau's results to all Serre curves further established that the average of the individual densities $C^{\Cyc}_{E,n,k}$ is equal to the average density $C^{\Cyc}_{n,k}$. This allows explicit determination of the values of $k$ that maximize and minimize $C^{\Cyc}_{E,n,k}$ for a fixed $n$ and computation of the corresponding extreme values.

Akbal and G\"{u}lo\v{g}lu also made the following interesting observation \cite[Proposition 1]{MR4504664}:
\begin{align}\label{AGobs}
    \text{There exists a prime } \ell \text{ such that } \Q(E[\ell]) \subset \Q(\zeta_n) \text{ and } \sigma_k|_{\Q(E[\ell])} \equiv 1 \hs \Longrightarrow \hs \lim_{x\to \infty}\pi^{\Cyc}_E(x;n,k) < \infty.
\end{align}
Thus, assuming GRH, if an elliptic curve $E$ satisfies $\Q(E[2]) \neq \Q$ and the condition in \eqref{AGobs}, the primes of cyclic reduction for $E$ exhibit a congruence class bias modulo $n$. Akbal and G\"ulo\v{g}lu also posed the question of whether the converse of \eqref{AGobs} holds. This was later shown to be false by Jones and the author through a concrete example:
\begin{Example}[\protect{\cite[Example 1.4]{jones2023acyclicity}}]\label{JLexa}
    Consider an elliptic curve given by
    $$E : Y^2 + XY + Y = X^3 + 32271697x - 1200056843302.$$
    This elliptic curve has conductor 71610, and its LMFDB label is $71610.s6$. We have
    \begin{itemize}
        \item $\Q(E[2]) \neq \Q$,
        \item $C_{E,8,k}^{\Cyc} = 0$ for $k \in \{3,5\}$,
        \item For each prime $\ell$, we have $\Q(E[\ell]) \not \subset \Q(\zeta_8)$.
    \end{itemize}
    \end{Example}
Furthermore, Jones and the author provided a complete characterization of when the density $C^{\Cyc}_{E,n,k}$ vanishes, thereby refining \eqref{AGobs} into a biconditional statement.
\begin{align*}
\begin{pmatrix}
    \text{There exists } d \in \N \text{ such that } \forall \sigma \in \Gal(\Q(E[d])\Q(\zeta_n)/\Q) \text{ with } \\
    \sigma |_{\Q(\zeta_n)} = \sigma_k, \exists \text{ a prime } \ell \mid d \text{ for which } \sigma|_{\Q(E[\ell])} = 1
\end{pmatrix} \iff C^{\Cyc}_{E,n,k} = 0.
\end{align*}
Building on this observation, they also constructed a modular curve that generates infinitely many isomorphism classes of elliptic curves, providing counterexamples to the question of Akbal and G\"ulo\v{g}lu. (See \cite[Theorem 1.7]{jones2023acyclicity}.)

\section{Prerequisites}\label{s3}
\subsection{Statistical Densities}\label{s3.1}
In this section, we express the average densities mentioned in \thmref{BSsummary} and \thmref{maintheorem} explicitly.

Define
$$C^{\Cyc} \coloneqq \prod_q \left( 1 - \frac{1}{q(q-1)(q^2-1)}\right) \approx 0.81375191\ldots ,$$
where the product runs over all primes $q$. By \thmref{BSsummary}, this quantity represents the average density of primes of cyclic reduction.

Let $\omega_r$ be a multiplicative function  defined at a prime power $q^j$ for $j \geq 1$ as follows:
\begin{equation}\label{omegardefine}
    \omega_r(q^j) \coloneqq \begin{cases} \displaystyle \frac{1}{q^{j-1}(q-1)} &\text{ if } r \not \equiv 1 \pmod {q^{\lceil j/2\rceil}} \\
\\
\displaystyle \frac{q^{\lfloor j/2 \rfloor + 1} + q^{\lfloor j/2\rfloor} - 1}{q^{j+\lfloor j/2 \rfloor -1}(q^2-1)} &\text{ if } r \equiv 1 \pmod {q^{\lceil j/2 \rceil}} \end{cases}.
\end{equation}
For any positive integer $m$, it follows directly from the definition that $\omega_r(m) \leq 1$. Given a positive integer $m$, we define
\begin{align}\label{mbar}
    \underline{m} \coloneqq \prod_{q^j \parallel m} q^{\lceil j/2\rceil}.
\end{align}

Observe that $\omega_r(q^j)$ depends only one whether $r$ is congruent to 1 modulo $q^{\lceil j/2 \rceil}$. Thus,
\begin{equation}\label{omegadepends}
    r \equiv s \mkern-12mu \pmod {\underline{m}}  \hs \Longrightarrow \hs \omega_r(m) = \omega_s(m).
\end{equation}
Define
$$C^{\Div(m)} \coloneqq \frac{1}{\phi(\underline{m})} \sum_{\substack{1 \leq r \leq \underline{m} \\ (r,\underline{m}) = 1}} \omega_r(m).$$
By \thmref{BSsummary}, this quantity represents the average density of primes of $m$-divisibility.

Next, define a multiplicative function $g$ at a prime power $q^j$ for $j\geq 1$ as 
\begin{equation}\label{definitionofg}
    g(q^j) \coloneqq \begin{cases} \displaystyle -\frac{1}{q(q^2-1)} &\text{ if } j = 1 \\ 0 & j \geq 2 \end{cases}
\end{equation}
Note that $q(q^2-1)$ is the size of $\SL_2(\Z/q\Z)$. Set
$$ C^{\Cyc}_{n,k} \coloneqq \frac{1}{\phi(n)}\prod_{\ell' \mid (n,k-1)}\left( 1 + g(\ell')\right) \prod_{\ell \nmid n}\left( 1 + \frac{g(\ell)}{\phi(\ell)}\right),$$
where $\ell'$ runs over all prime factors of $\gcd(n,k-1)$ and $\ell$ runs over all primes coprime to $n$. The infinite product converges absolutely, since
$$\left|\frac{g(\ell)}{\phi(\ell)}\right| = \frac{1}{|\GL_2(\Z/\ell \Z)|} \sim \frac{1}{\ell^4}, \hs \text{ as }\hs \ell \to \infty.$$
By \thmref{maintheorem}, this quantity represents the average density of primes $p \equiv k \pmod n$ of cyclic reduction.

Finally, for positive integers $n$ and $m$, set
\begin{equation}\label{definitionofM}
    \mathcal{M} \coloneqq \lcm(n,\underline{m}).
\end{equation}
Define
$$C^{\Div(m)}_{n,k} \coloneqq \frac{1}{\phi(\mathcal{M})} \sum_{\substack{1 \leq r \leq \mathcal{M} \\ (r, \mathcal{M}) = 1 \\ r \equiv k (n)}}\omega_r(m).$$
By \thmref{maintheorem}, this quantity represents the average density of primes $p \equiv k \pmod n$ of $m$-divisibility.

\subsection{Modifying Banks and Shparlinski's methods}\label{s3.2}
In this section, we modify the main argument of \cite{MR2570668} to compute the asymptotic formulas for
$$\frac{1}{4AB}\sum_{|a| \leq A} \sum_{|b| \leq B} \pi_{a,b}^{\Cyc}(x;n,k) \hs \text{ and } \hs  \frac{1}{4AB}\sum_{|a| \leq A} \sum_{|b| \leq B} \pi^{\Div(m)}_{a,b}(x;n,k).$$
Take $\epsilon > 0$ and $x > 0$. Set $A \coloneqq A(x)$ and $B \coloneqq B(x)$ as positive integers satisfying
$$x^{\epsilon} \leq A, B \leq x^{1-\epsilon}, \hs \hs x^{1+\epsilon} \leq AB.$$
This condition becomes empty if $\epsilon > 1/3$.

For simplicity, let $\mathcal{X}$ denote either $\Cyc$ or $\Div(m)$. Define
    $$N^{\mathcal{X}}_{n,k}(A,B;x) \coloneqq \sum_{|a| \leq A} \sum_{|b| \leq B} \pi^{\mathcal{X}}_{a,b}(x;n,k).$$
Fix a prime $p > 3$. Any elliptic curve over $\F_p$ can be written as
$$\mathcal{E}^{\alpha,\beta} : Y^2 \equiv X^3 + \alpha X + \beta,$$
for some $\alpha,\beta \in \F_p$ with $4\alpha^3+27\beta^2 \neq 0$. Define
\begin{align*}
    \Cyc_p &\coloneqq \{(\alpha,\beta) \in \F_p^\times \times \F_p^\times : 4\alpha^3 + 27\beta^2 \neq 0 \text{ and } \mathcal{E}^{\alpha,\beta}(\F_p) \text{ is cyclic}\}, \\
    \Div_p(m) & \coloneqq \{(\alpha,\beta) \in \F_p^\times \times \F_p^\times : 4\alpha^3 + 27\beta^2 \neq 0 \text{ and } m \text{ divides } |\mathcal{E}^{\alpha,\beta}(\F_p)|\}.
\end{align*}
Likewise, let $\mathcal{X}_p$ denote either $\Cyc_p$ or $\Div_p(m)$. We define
\begin{align*}
    M_p(\mathcal{X}_p;A,B) \coloneqq \#\{(a,b) \in \Z \times \Z : E^{a,b} \in \mathcal{F}(A,B) \text{ and }  (\bar{a}, \bar{b}) \in \mathcal{X}_p\},
\end{align*}
where $\bar{a}$ and $\bar{b}$ denote $a \pmod p$ and $b \pmod p$, respectively.

\begin{Remark}\label{observation}
Let $\mathcal{F}_p$ denote the family of all models $\mathcal{E}^{\alpha,\beta}$ of elliptic curves over $\F_p$ and $\mathcal{F}_p^{\mathcal{X}}$, the subfamily satisfying the property $\mathcal{X}$. Note that

    
    $$N^{\mathcal{X}}_{n,k}(A,B;x) = \sum_{E \in \mathcal{F}} \sum_{\substack{p \leq x \\ p \equiv k (n) \\ \widetilde{E}_p \in \mathcal{F}_p^{\mathcal{X}}}}1.$$
    Let $M = \max\{A,B\}$. Since $M \leq x^{1-\epsilon}$,  we have
    $$N^{\mathcal{X}}_{n,k}(A,B;x) = \sum_{E \in \mathcal{F}} \left( \sum_{\substack{2M + 1 < p \leq x \\ p \equiv k (n) \\ \widetilde{E}_p \in \mathcal{F}_p^{\mathcal{X}}}} 1 + \sum_{\substack{p \leq 2M + 1 \\ p \equiv k (n) \\ \widetilde{E}_p \in \mathcal{F}_p^{\mathcal{X}}}} 1\right).$$
    For the cyclicity case, we take $n \leq \log^K x$, while for the divisibility case, we take $n$ and $m$ such that $\mathcal{M} \leq \log^{K+1} x$, where $\mathcal{M}$ is  defined in \eqref{definitionofM}. By definition, $n \leq \mathcal{M}$ which implies $n \leq \log^{K+1} x$. Moreover, with $x^\epsilon \leq M$, it follows that $n \ll_{\epsilon, K} M$ in both cases. Since $|\mathcal{F}| \leq 4AB$, we obtain
    $$\sum_{E \in \mathcal{F}} \sum_{\substack{p \leq 2M + 1 \\ p \equiv k (n) \\ \widetilde{E}_p \in \mathcal{F}_p^{\mathcal{X}}}} 1 \leq \sum_{E \in \mathcal{F}} \left(\frac{2M+1}{n} + 1 \right) \leq \frac{8MAB+4AB}{n} + 4AB \ll_{\epsilon,K} \frac{MAB}{n}.$$
   For $p > 2M+1$, the reduction map $\mathcal{F} \to \mathcal{F}_p$ is injective. Suppose $\Eab \in \mathcal{F}$ with $\Eabp \in \mathcal{F}_p^{\mathcal{X}}$. Then, either $\bar{a}\bar{b} \in \F_p^\times$ or $\bar{a}\bar{b} 
 = 0$. The former case corresponds to $(\bar{a},\bar{b}) \in \mathcal{X}_p$, while the latter occurs when either $a = 0$ or $b = 0$. Since there are at most $4M$ such elliptic curves, we get
    $$\sum_{E \in \mathcal{F}} \sum_{\substack{2M + 1 < p \leq x \\ p \equiv k (n) \\ \widetilde{E}_p \in \mathcal{F}_p^{\mathcal{X}}}} 1 \leq \sum_{\substack{2M + 1 < p \leq x \\ p \equiv k(n)}} \left( \sum_{\substack{|a| \leq A \\ |b| \leq B \\ (\bar{a},\bar{b}) \in \mathcal{X}_p}} 1 +4M\right) \leq \sum_{\substack{2M + 1 < p \leq x \\ p \equiv k (n)}} M_p (\mathcal{X}_p;A,B) + \frac{4Mx}{n} + 4M.$$
    Since $AB \geq x^{1+\epsilon}$, it follows that $Mx \leq MAB$. Thus, we obtain
    \begin{equation}\label{sumfinal}
        N^{\mathcal{X}}_{n,k}(A,B;x) = \sum_{\substack{2M + 1 < p \leq x \\ p \equiv k (n)}} M_p(\mathcal{X}_p; A,B) + O_{\epsilon,K}\left(\frac{MAB}{n}\right).
    \end{equation}
\end{Remark}
Banks and Shparlinski determined the size of $M_p(\mathcal{X}_p;A,B)$ in terms of $p$ with a manageable error bound. Let $\mathfrak{X}_p$ be the set of all multiplicative characters of $\F_p$, and let $\chi_0$ be the principal character modulo $p$. Define $d_p \coloneqq \gcd(p-1,6)$ and $e_p \coloneqq \gcd(p-1,4)$. For a positive integer $N$, set
\begin{align*}
    \sigma_p(N) \coloneqq \max_{\substack{\chi \in \mathfrak{X}_p\setminus\{\chi_0\} \\ \chi^{d_p} = \chi_0}} \left\{ 1 , \left| \sum_{n=1}^N \chi(n) \right| \right\}, \hs \rho_p(N) \coloneqq \max_{\substack{\chi \in \mathfrak{X}_p\setminus\{\chi_0\} \\ \chi^{e_p} = \chi_0}} \left\{ 1 , \left| \sum_{n=1}^N \chi(n)\right|\right\}
\end{align*}
Define
\begin{equation}\label{EABp}
    \mathcal{E}(A,B;p) \coloneqq \min\{A\sigma_p(B) + A^{1/2}Bp^{1/4} + A^{1/2}B^{1/2}p^{1/2}, B\rho_p(A) + AB^{1/2}p^{1/4} + A^{1/2}B^{1/2}p^{1/2}\}
\end{equation}

\begin{Theorem}[\protect{\cite[Theorem 13]{MR2570668}}]\label{BSPrelim}
Let $p > 3$ be a prime and $A,B$ be integers satisfying $1 \leq A,B \leq (p-1)/2$. Let $\mathcal{S}_p \subseteq \F_p \times \F_p$ be a set closed under isomorphism, i.e., if $(a,b) \in \mathcal{S}_p$ and $\mathcal{E}^{a,b}(\F_p) \cong \mathcal{E}^{c,d}(\F_p)$, then $(c,d) \in \mathcal{S}_p$. Then, for any choice of $\epsilon' > 0$, the following bound holds uniformly.
$$M_p(\mathcal{S}_p; A,B) - \frac{4AB}{p^2}\#\mathcal{S}_p \ll_{\epsilon'} \mathcal{E}(A,B;p)p^{\epsilon'}, \hs \text{ as } p \to \infty$$
\end{Theorem}

Since $\Cyc_p$ and $\Div_p(m)$ are closed under isomorphism, it suffices to compute their cardinalities with manageable error bounds, a task accomplished by Vl\v{a}du\c{t} and Howe, respectively.

Define
\begin{align}\label{definitionofvartheta}
    \vartheta_p \coloneqq \prod_{q \mid p-1} \left( 1 - \frac{1}{q(q^2-1)}\right),
\end{align}
where the product is taken over all primes $q$ dividing $p-1$.
\begin{Lemma}[\protect{\cite[Theorem 6.1]{MR1667099}}]\label{cardinalityofCp}
For all primes $p$ and $\epsilon' > 0$,
$$\left| \# \Cyc_p-\vartheta_p p^2 \right| \ll_{\epsilon'} p^{3/2+\epsilon'},\hs \text{as } p \to \infty.$$
\end{Lemma}
\begin{Lemma}[\protect{\cite[Theorem 1.1]{MR1204781}}]
\label{cardinalityofDp}
For all primes $p$, positive integers $m$, and $\epsilon' > 0$,
$$\left| \#\Div_p(m) - \omega_p(m)p^2 \right| \ll_{\epsilon'} m^{1+\epsilon'}p^{3/2}, \hs \text{ as } m \to \infty.$$
\end{Lemma}
(See also Lemma 10 and Lemma 11 in \cite{MR2570668}.)

Let us check the cyclicity case first. Applying \lemref{cardinalityofCp} to \thmref{BSPrelim}, we obtain
\begin{align}\label{MpCp1}
    |M_p(\Cyc_p;A,B) - 4AB \vartheta_p | \ll_{\epsilon'} 4ABp^{-1/2+\epsilon'} + \mathcal{E}(A,B;p)p^{\epsilon'}.
\end{align}
Since $p > 2M+1$, we have $AB < p^2$, and hence $ABp^{-1/2} < A^{1/2}B^{1/2}p^{1/2}$. Thus, the term $4ABp^{-1/2 + \epsilon'}$ can be absorbed into the error term $\mathcal{E}(A,B;p)p^{\epsilon'}$. Summing \eqref{MpCp1} over $p \leq x$ for which $p \equiv k \pmod n$, we obtain
\begin{align*}
    \left|\sum_{\substack{ p\leq x \\ p \equiv k (n)}} M_p(\Cyc_p;A,B) - 4AB\sum_{\substack{p \leq x \\ p \equiv k (n)}} \vartheta_p \right| \ll_{\epsilon'} \sum_{\substack{p \leq x \\ p \equiv k (n)}} \mathcal{E}(A,B;p)p^{\epsilon'} \leq x^{\epsilon'} \sum_{\substack{p \leq x \\ p \equiv k (n)}} \mathcal{E}(A,B;p).
\end{align*}
Using \eqref{sumfinal}, we have
\begin{align}\label{MpCp3}
    \left| N_{n,k}^{\Cyc}(A,B;x) - 4AB\sum_{\substack{p \leq x \\ p \equiv k (n)}} \vartheta_p \right| \ll_{\epsilon'} x^{\epsilon'} \sum_{\substack{p \leq x \\ p \equiv k (n)}} \mathcal{E}(A,B;p) + O_{\epsilon,K} \left( \frac{MAB}{n}\right).
\end{align}

Now, let us check the divisibility case. Applying \lemref{cardinalityofDp} to \thmref{BSPrelim}, we obtain
\begin{align}\label{MpDp1}
    |M_p(\Div_p(m);A,B) - 4AB\omega_p(m)| \ll_{\epsilon'} 4ABm^{1+\epsilon'}p^{-1/2} + \mathcal{E}(A,B;p)p^{\epsilon'}.
\end{align}
By \eqref{mbar} and \eqref{definitionofM}, it follows that $m \leq \underline{m}^2 \leq \mathcal{M}^2 \leq \log^{2K+2} x$. Thus, we have 
\begin{equation}\label{auxil}
   m \leq \log^{2K+2}x \ll_{\epsilon',K} x^{\epsilon'} \implies 4ABm^{1+\epsilon'}p^{-1/2} \ll_{\epsilon',K} x^{\epsilon'} \mathcal{E}(A,B;p).
\end{equation}
Therefore, \eqref{sumfinal}, \eqref{MpDp1}, and \eqref{auxil} give
\begin{align}\label{MpDp2}
    \left|N^{\Div(m)}_{n,k}(A,B;x) - 4AB \sum_{\substack{p \leq x \\ p \equiv k (n)}} \omega_p(m)\right| \ll_{\epsilon',K} x^{\epsilon'} \sum_{\substack{p \leq x \\ p \equiv k (n)}} \mathcal{E}(A,B;p) + O_{\epsilon,K}\left( \frac{MAB}{n}\right).
\end{align}

\subsection{On the sum of \texorpdfstring{$\mathcal{E}(A,B;p)$}{EABp}}\label{s3.3}
In this section, we examine the sum of $\mathcal{E}(A,B;p)$ and refine the results of Banks and Shparlinski. For notational convenience, we define
\begin{align*}
    \alpha_p & \coloneqq  A \sigma_p(B) + A^{1/2}Bp^{1/4} + A^{1/2}B^{1/2}p^{1/2}, \\
 \beta_p & \coloneqq B\rho_p(A) + AB^{1/2}p^{1/4} + A^{1/2}B^{1/2}p^{1/2}.
\end{align*}
By \eqref{EABp}, the definition of $\mathcal{E}(A,B;p)$, we obtain
\begin{align*}
    \sum_{\substack{p \leq x \\ p \equiv k (n)}} \mathcal{E}(A,B;p) = \sum_{\substack{p \leq x \\ p \equiv k (n)}} \min\{\alpha_p, \beta_p\} \leq \min \left\{\sum_{\substack{p \leq x \\ p \equiv k (n)}} \alpha_p, \sum_{\substack{p \leq x \\ p \equiv k (n)}} \beta_p \right\}.
\end{align*}
We apply the following bound:
\begin{equation}\label{upperboundofEABp}
    \sum_{\substack{p \leq x \\ p \equiv k (n)}}\mathcal{E}(A,B;p) \leq \begin{cases}
    \displaystyle \sum_{\substack{p \leq x \\ p \equiv k (n)}}\alpha_p & \text{ if } A \geq B, \\ \\
    \displaystyle \sum_{\substack{p \leq x \\ p \equiv k(n)}}\beta_p & \text{ if }A < B.
\end{cases}
\end{equation}
To bound $\alpha_p$ and $\beta_p$, we use results by M. Z. Garaev on upper bounds for character sums. (See also \cite[Lemma 3]{MR2570668}.)
\begin{Theorem}[\protect{\cite[Theorem 10]{MR2235360}}] \label{GaraevTheorem}
Fix $\epsilon > 0$ and $\eta > 0$. If $x$ is sufficiently large, then for all $N \geq x^{\epsilon}$, all primes $p \leq x$ with at most $O_\epsilon(x^{3/4 + 4\eta + o_{\epsilon}(1)})$ exceptions as $x \to \infty$, and all characters $\chi \in \mathfrak{X}_p\setminus \{\chi_0\}$, we have
$$\left| \sum_{n = 1}^N \chi(n)\right| \leq N^{1-\eta}.$$
\end{Theorem}
    Banks and Shparlinski assumed that the little $o$-term depends on both $\epsilon$ and $\eta$. However, upon reviewing Garaev's proof of the theorem, one could confirm that the little $o$-term depends only on $\epsilon$. This confirmation is crucial because, if the little $o$-term depended on both $\epsilon$ and $\eta$, it could impose restrictions on their choices. Since $\epsilon$ and $\eta$ are chosen independently in our proof of the main theorem, ensuring that the little $o$-term depends only on $\epsilon$ guarantees that the error remains controlled.

Fix a prime $p$ and $\eta > 0$. (Recall that \thmref{GaraevTheorem} provides a meaningful result only if $\eta < 1/16$.) The theorem implies that for any $p \leq x$ with $p \equiv k \pmod n$, we have
$$\sigma_p(B) \leq B^{1-\eta},$$
with at most $O_{\epsilon}(x^{3/4 + 4\eta + o_{\epsilon}(1)})$ exceptions. (If $p$ is one of these exceptional primes, we use the trivial bound $\sigma_p(B) \leq B$.) Thus, we obtain
$$A\sum_{\substack{p \leq x \\ p \equiv k (n)}} \sigma_p(B) \ll_\epsilon ABx^{3/4+4\eta+o_{\epsilon}(1)} + \frac{1}{n} AB^{1-\eta}x.$$
Assume $A \geq B$. Using \eqref{upperboundofEABp}, we get
\begin{equation}\label{ifA>B}
    x^{\epsilon'}\sum_{\substack{p \leq x \\ p \equiv k (n)}} \mathcal{E}(A,B;p) \ll_\epsilon ABx^{3/4+4\eta+\epsilon'+o_{\epsilon}(1)} + \frac{x^{\epsilon'}}{n} \left( AB^{1-\eta} x + A^{1/2}Bx^{5/4} + A^{1/2}B^{1/2}x^{3/2}\right).
\end{equation}
Otherwise, for $B > A$, we get
$$x^{\epsilon'}\sum_{\substack{p \leq x \\ p \equiv k (n)}} \mathcal{E}(A,B;p) \ll_\epsilon ABx^{3/4+4\eta+\epsilon'+o_{\epsilon}(1)} + \frac{x^{\epsilon'}}{n} \left( A^{1-\eta}B x + AB^{1/2}x^{5/4} + A^{1/2}B^{1/2}x^{3/2}\right).$$
\begin{Proposition}\label{errorn}
    Set $\epsilon > 0$ and $x > 0$. Let $A \coloneqq A(x)$ and $B \coloneqq B(x)$ satisfy \eqref{conditionsonAB}. Choose $\sigma$ and $\eta$ such that
    $$ 0 < \frac{\sigma}{\epsilon} < \eta < \frac{1}{16} \hs \text{ and } \hs \sigma + 4\eta < \frac{1}{4}.$$
    Let $\epsilon'$ and $\delta$ satisfy
    $$0 < \epsilon' < \min \left\{\epsilon\eta - \sigma, \frac{3\epsilon}{16}\right\} \hs \text{ and } \hs 0 < \delta < \frac{1}{4} - \sigma - 4\eta.$$
    Then, for any $n \leq x^{1/4-\sigma-4\eta-\delta}$, we have
        $$x^{\epsilon'}\sum_{\substack{p \leq x \\ p \equiv k (n)}}\mathcal{E}(A,B;p) \ll_{\epsilon,\delta} \frac{1}{n}ABx^{1-\sigma},$$
    for sufficiently large $x$.
\end{Proposition}
\begin{proof}
Suppose $A \geq B$. We will verify that each term on the right-hand side of \eqref{ifA>B} can be bounded by $ABx^{1-\sigma}/n$, for sufficiently large $x$.

Given the conditions, we have $\sigma < \epsilon/16$ and $\epsilon' < 3\epsilon/16$. Therefore, $\sigma + \epsilon' < \epsilon/4$, so $x^{\sigma + \epsilon'} \leq x^{\epsilon/4}$ for $x \geq 1$. Since $x^{1+\epsilon} \leq AB$, we have
$$x^{1/2 + \sigma + \epsilon'} \leq x^{1/2+\epsilon/2} \leq A^{1/2}B^{1/2} \hs \implies \hs \frac{1}{n}A^{1/2}B^{1/2}x^{3/2+\epsilon'} \leq \frac{1}{n}ABx^{1-\sigma},$$
for $x \geq 1$.

Given $A \geq B$ and $AB \geq x^{1+\epsilon}$, we have $A \geq x^{1/2 + \epsilon/2}$. Thus, we have $$x^{1/4 + \sigma + \epsilon'} \leq x^{1/4 + \epsilon/4} \leq A^{1/2} \hs \Longrightarrow \hs \frac{1}{n}A^{1/2}Bx^{5/4+\epsilon'} \leq \frac{1}{n}ABx^{1-\sigma},$$
for $x \geq 1$.

Since $\epsilon' + \sigma < \epsilon \eta$ and $x^\epsilon \leq B$, we have 
$$x^{\epsilon' + \sigma} \leq x^{\epsilon \eta} \leq B^\eta \implies \frac{1}{n} AB^{1-\eta} x^{1+\epsilon'} \hs \leq \hs \frac{1}{n} ABx^{1-\sigma},$$
for $x \geq 1$. Therefore, for $A \geq B$, we have
$$\frac{x^{\epsilon'}}{n}\left(AB^{1-\eta}x + A^{1/2}Bx^{5/4} + A^{1/2}B^{1/2}x^{3/2}\right) \ll \frac{1}{n} ABx^{1-\sigma}.$$

If $B > A$, we can argue similarly to obtain
$$\frac{x^{\epsilon'}}{n}\left(A^{1-\eta}Bx + AB^{1/2}x^{5/4} + A^{1/2}B^{1/2}x^{3/2} \right) \ll \frac{1}{n}ABx^{1-\sigma}.$$

Since $\delta > 0$, we have $x^{o_{\epsilon}(1)} \ll_{\epsilon,\delta} x^\delta$. Given the condition on $n$, we see that
\begin{align*}
   nx^{3/4 + 4\eta + o_\epsilon(1)} \leq x^{1 - \sigma -\delta + o_{\epsilon}(1)} \ll_{\epsilon,\delta} x^{1-\sigma} \hs \implies \hs ABx^{3/4 + 4\eta + o_{\epsilon}(1)} \ll_{\epsilon,\delta} \frac{1}{n}ABx^{1-\sigma},
\end{align*}
for sufficiently large $x$.
\end{proof}

\begin{Proposition}\label{errorn1}
       Set $\epsilon > 0$ and $x > 0$. Let $A \coloneqq A(x)$ and $B \coloneqq B(x)$ satisfy \eqref{conditionsonAB}. Choose $\sigma$ and $\eta$ such that
    $$ 0 < \frac{\sigma}{\epsilon} < \eta < \frac{1}{16} \hs 
    \text{ and } \hs \sigma + 4\eta < \frac{1}{4}.$$
    Let $\epsilon'$ satisfy
    $$0 < \epsilon' < \min \left\{\epsilon\eta - \sigma, \frac{3\epsilon}{16}\right\}.$$ 
    Take $K > 0$. For any $n \leq \log^{K} x$, we have
    $$x^{\epsilon'} \sum_{\substack{p \leq x \\ p \equiv k (n)}}\mathcal{E}(A,B;p) \ll_{\epsilon, \sigma, \eta, K} \frac{1}{n}ABx^{1-\sigma},$$
    for sufficiently large $x$.
\end{Proposition}
\begin{proof}
    The proof follows similarly to the previous case. We only need to establish an upper bound for $ABx^{3/4+4\eta + o_{\epsilon}(1)}$. Since $n \leq \log^K x$ and $1/4 - \sigma - 4\eta > 0$, we have
    $$nx^{o_{\epsilon}(1)} \leq x^{o_{\epsilon}(1)} \log^K x \ll_{\epsilon, \sigma, \eta, K} x^{1/4 - \sigma - 4\eta} \hs \implies \hs ABx^{3/4 + 4\eta + o_{\epsilon}(1)} \ll_{\epsilon,\sigma,\eta,K} \frac{1}{n}ABx^{1-\sigma},$$
    for sufficiently large $x$.
\end{proof}

\begin{Corollary}\label{errornsmall}
     Set $\epsilon > 0$ and $x > 0$. Let $A \coloneqq A(x)$ and $B \coloneqq B(x)$ satisfy \eqref{conditionsonAB}. Set $\sigma < \epsilon/(16+4\epsilon)$, $K > 0$, and $n \leq \log^K x$. Then for sufficiently small $\epsilon' > 0$, we have
    $$x^{\epsilon'} \sum_{\substack{p \leq x \\ p \equiv k (n)}} \mathcal{E}(A,B;p) \ll_{\epsilon, \sigma, K} \frac{1}{n}ABx^{1-\sigma},$$
    for sufficiently large $x$.
\end{Corollary}
\begin{proof}
    The conditions in Proposition \ref{errorn1} provide upper bounds on $\sigma$:
    $$\frac{\sigma}{\epsilon} < \eta \hs \text{ and } \hs \sigma + 4\eta < \frac{1}{4} \hs \implies \hs \sigma < \min\left\{\epsilon\eta, \frac{1}{4} - 4\eta\right\}.$$
    Thus, the largest upper bound on $\sigma$ can be attained if those two quantities are equal:
    $$\epsilon \eta = \frac{1}{4} - 4\eta \hs \implies \hs \eta = \frac{1}{16+4\epsilon}.$$
    Thus, by choosing $\sigma < \epsilon/(16+4\epsilon)$, we ensure that $\sigma$ and $\eta$ satisfy the conditions in Proposition \ref{errorn1}. Set $\epsilon'= \frac{1}{2}\min\{\epsilon \eta - \sigma, 3\epsilon/16\}$. By Proposition \ref{errorn1}, we obtain the desired results.
\end{proof}
Recall that $\epsilon < 1/3$, due to \eqref{conditionsonAB}. Therefore, the supremum of $\sigma$ is $1/52$.

\begin{Corollary}\label{errornlarge}
    Set $\epsilon > 0$ and $x > 0$. Let $A \coloneqq A(x)$ and $B \coloneqq B(x)$ satisfy \eqref{conditionsonAB}. Set $\tau \leq 1/4 -11\epsilon/16$. For any $K > 0$ and $n \leq x^\tau$, we have
    $$x^{\epsilon'} \sum_{\substack{p \leq x \\ p \equiv k (n)}} \mathcal{E}(A,B;p) \ll_{\epsilon, \tau, K} \frac{ABx}{n\log^{K+1} x},$$
    for sufficiently large $x$.
\end{Corollary}
\begin{proof}
    We define the parameters $\eta, \sigma$, $\delta$, and $\epsilon'$ as follows:
\begin{align*}
    \eta = \frac{\epsilon}{8}, \hs \hs \sigma = \frac{\epsilon^2}{16}, \hs \hs \delta = \frac{\epsilon}{8}, \hs \hs \epsilon' = \frac{\epsilon^2}{32}.
\end{align*}
It is straightforward to check that $\sigma$ and $\eta$ meet the required conditions in \propref{errorn}. Since $\epsilon < 1/3$, we have
$$\epsilon' < \epsilon \eta - \sigma=  \frac{\epsilon^2}{16}  < \frac{3\epsilon}{16} \hs \text{ and } \hs \sigma + 4\eta + \delta = \frac{\epsilon^2}{16} + \frac{\epsilon}{2} + \frac{\epsilon}{8} < \frac{11\epsilon}{16} < \frac{11}{48}< \frac{1}{4}.$$
Thus, $\epsilon'$ and $\delta$ also meet the conditions. Hence, we have
$$\tau \leq \frac{1}{4} - \frac{11\epsilon}{16}< \frac{1}{4} - \sigma - 4\eta - \delta \hs \implies \hs x^\tau \leq x^{1/4-\sigma-4\eta-\delta},$$
for any $x \geq 1$. Therefore, for any choice of $n \leq x^\tau$, by Proposition \ref{errorn}, we obtain
      $$x^{\epsilon'} \sum_{\substack{p \leq x \\ p \equiv k(n)}}\mathcal{E}(A,B;p) \ll_{\epsilon, \tau} \frac{1}{n}ABx^{1-\sigma},$$
      for sufficiently large $x$. Finally, observe that
      $$\log^{K+1} x \ll_{\epsilon, K} x^\sigma \implies \frac{1}{n}ABx^{1-\sigma} \ll_{\epsilon,K} \frac{ABx}{n\log^{K+1}x} \implies x^{\epsilon'} \sum_{\substack{p \leq x \\ p \equiv k(n)}}\mathcal{E}(A,B;p) \ll_{\epsilon, \tau, K} \frac{ABx}{n\log^{K+1} x}, $$
      for sufficiently large $x$.
\end{proof}

Recall the extra $O_{\epsilon, K}(MAB/n)$ term in \eqref{MpCp3} and \eqref{MpDp2}, where $M = \max\{A,B\} \leq x^{1-\epsilon}$. As we choose $\sigma < \epsilon/(16+4\epsilon)$, we have
\begin{equation}\label{Mterm}
    M \leq x^{1-\sigma} \implies \frac{MAB}{n} \leq \frac{1}{n}ABx^{1-\sigma} \hs \text{ and } \hs  M \log^{K+1} x \ll_{\epsilon, K} x \implies \frac{MAB}{n} \ll_{\epsilon,K} \frac{ABx}{n\log^{K+1} x},
\end{equation}
for sufficiently large $x$. Let $\eta$ and $\epsilon'$ be defined as in the proof of Corollary \ref{errornsmall}. Using \eqref{MpCp3}, \eqref{Mterm}, and Corollary \ref{errornsmall}, we deduce the following lemma.
\begin{Lemma}\label{NCnk0}
    Let $\epsilon > 0$ and $x > 0$. Let $A \coloneqq A(x)$ and $B \coloneqq B(x)$ be positive integers satisfying
    $$x^{\epsilon} \leq A, B \leq x^{1-\epsilon}, \hs AB \geq x^{1+\epsilon}.$$
    Fix $\sigma < \epsilon/(16+4\epsilon)$ and $K > 0$. Then, for any $n \leq \log^K x$, we have
    $$N^{\Cyc}_{n,k}(A,B;x) = 4AB\sum_{\substack{p \leq x \\ p \equiv k (n)}} \vartheta_p + O_{\epsilon, \sigma, K} \left(\frac{1}{n}ABx^{1-\sigma} \right),$$
    for sufficiently large $x$.
\end{Lemma}
Similarly, let $\eta, \sigma, \delta, \epsilon'$ be set as in the proof of Corollary \ref{errornlarge}. Using \eqref{MpDp2}, \eqref{Mterm}, and Corollary \ref{errornlarge}, we obtain the following lemma.

\begin{Lemma}\label{NDnk0}
    Let $\epsilon > 0$ and $x > 0$. Let $A \coloneqq A(x)$ and $B \coloneqq B(x)$ be positive integers satisfying
    $$x^{\epsilon} \leq A, B \leq x^{1-\epsilon}, \hs AB \geq x^{1+\epsilon}.$$
    Set $\tau \leq 1/4 - 11\epsilon/16$. For any $K > 0$ and $n \leq x^\tau$, we have
    $$N^{\Div(m)}_{n,k}(A,B;x) = 4AB\sum_{\substack{p \leq x \\ p \equiv k (n)}} \omega_p(m) + O_{\epsilon, \tau, K}\left( \frac{ABx}{n\log^{K+1}x}\right),$$
    for sufficiently large $x$.
\end{Lemma}
Thus, we are reduced to finding the asymptotic formulas for
$$\sum_{\substack{p \leq x \\ p \equiv k (n)}} \vartheta_p \hs \text{ and } \hs \sum_{\substack{p \leq x \\ p \equiv k (n)}}\omega_p(m).$$
\begin{Remark}
   There are two different strategies for choosing $\epsilon$, each offering distinct advantages and drawbacks. One approach is to maximize $\epsilon$, which  sharpens the error bound due to a larger upper bound of $\sigma$ as seen in Corollary \ref{errornsmall}. However, this comes at the cost of restricting the values of $A$ and $B$, pushing them closer together. In the extreme case, setting $\epsilon = 1/3$ results in $A = B$, which reduces the two-parameter family $\mathcal{F}(A,B)$ to a single-parameter family. 

   On the other hand, one could choose $\epsilon$ to be as small as possible, which allows more flexibility in selecting $A$ and $B$ and permits larger values for $n$ as seen in Corollary \ref{errornlarge}. However, this comes with a worse error bound.

    In proving the asymptotic behaviors of the sums involving $\vartheta_p$ and $\omega_p(m)$, we leverage different theorems. For the former, the Bombieri-Vingoradov theorem, with its smaller error bound, allows
us to maximize $\epsilon$ and hence so thus $\sigma$. In contrast, for the latter, we rely on the Siegel-Walfisz theorem, which features a comparatively large error term, necessitating a choice of smaller $\epsilon$.
\end{Remark}

The following lemmas describe the asymptotic formulas for the sums of $\vartheta_p$ and $\omega_p(m)$, which will be proven in the next chapter.
\begin{Lemma}\label{sumofvarthetap}
Let $g$ be the function defined in \eqref{definitionofg} and $K > 0$. Let $n \leq \log^K x$ be a positive integer and $k$ be coprime to $n$. Then, for any $\delta > 0$, we have
$$\sum_{\substack{p \leq x \\ p \equiv k (n)}} \vartheta_p = \frac{\li (x)}{\phi (n)} \prod_{\ell'\mid (n,k-1)} (1+g(\ell'))\prod_{\ell \nmid n} \left( 1 + \frac{g(\ell)}{\phi(\ell)}\right) + O_{\delta,K}\left(x^{1/2 + \delta}\right),$$
for sufficiently large $x$.
\end{Lemma}
\begin{Lemma}\label{sumofomegap}
Let $n$ and $m$ be positive integers and $k$ be coprime to $n$. Let $\mathcal{M}$ be as defined in \eqref{definitionofM}. For any $K > 0$ with $\mathcal{M} \leq \log^{K+1} x$, there exists a constant $c$ that only depends on $K$ such that
$$\sum_{\substack{p \leq x \\ p \equiv k(n)}} \omega_p(m) = \frac{\li (x)}{\phi(\mathcal{M})} \sum_{\substack{ 1 \leq r \leq \mathcal{M} \\ (r,\mathcal{M}) =1 \\ r \equiv k (n)}} \omega_r(m) + O \left(m x \exp \left( -c \sqrt{\log x}\right) \right),$$
for sufficiently large $x$.
\end{Lemma}

We conclude this chapter with the proof of the average Sato-Tate conjecture, assuming primes lie in an arithmetic proression.
\begin{Remark}\label{STACC}
    For $0 \leq \alpha < \beta \leq \pi$, let
    \begin{align*}
        \ST_p(\alpha,\beta) &\coloneqq \left\{(a,b) \in \F_p^\times \times \F_p^\times : \alpha \leq \psi_{\mathcal{E}^{a,b}} \leq \beta \right\},
    \end{align*}
    where we define $\psi_{\mathcal{E}^{a,b}} \in [0,\pi]$ satisfying 
    $$p+1-\#\mathcal{E}^{a,b}(\F_p) = 2\sqrt{p} \cos \psi_{\mathcal{E}^{a,b}}.$$
    Banks and Shparlinski \cite[Lemma 9]{MR2570668} proved that
    $$ \max_{0 \leq \alpha < \beta \leq \pi} \left|\#\ST_p(\alpha,\beta) - \mu_{\ST}(\alpha,\beta)p^2\right| \ll p^{7/4},$$
    where $\mu_{\ST}$ is as defined in \eqref{SatoTatemeasure}. Using \thmref{BSPrelim} and the argument sketched in \sectionref{s3.2}, one could derive
    \begin{equation}\label{Keyidea}
        \left| M_p(\ST_p(\alpha,\beta); A,B) - 4AB\mu_{\ST}(\alpha,\beta) \right| \ll_{\epsilon'} \mathcal{E}(A,B;p)p^{\epsilon'},
    \end{equation}
    Define
    $$\pi_{a,b}^{\ST(\alpha,\beta)}(x;n,k) \coloneqq \#\left\{ p \leq x : p \nmid N_{\Eab} , \psi_{\Eab}(p) \in [\alpha,\beta], p \equiv k \neghs \pmod n \right\},$$
    where $\psi_{\Eab}(p)$ is as defined in \eqref{Frobangle}. Summing \eqref{Keyidea} over $p \leq x$ and $p \equiv k \pmod n$ and using \corref{errornsmall}, we obtain the following result: Let $x > 0, \epsilon > 0$, and $K > 0$, and let $A$ and $B$ satisfy \eqref{conditionsonAB}. Set $n \leq \log^K x$ and $\sigma < \epsilon/(16+4\epsilon)$. Then, we have
    $$\frac{1}{4AB} \sum_{|a| \leq A} \sum_{|b| \leq B} \pi^{\ST(\alpha,\beta)}_{\Eab} (x;n,k) = \mu_{\ST}(\alpha,\beta)\pi(x;n,k) + O_{\epsilon, \sigma, K}\left( \frac{x^{1-\sigma}}{n}\right),$$
    for sufficiently large $x$. In particular, there is no average congruence class bias for primes whose Frobenius angle lies in a closed interval $[\alpha,\beta]$.
\end{Remark}

\section{Proofs of \lemref{sumofvarthetap} and \lemref{sumofomegap}}\label{s4}
\subsection{Mean behavior of an Arithmetic Function}\label{s4.1} We begin by exploring the mean behavior of an arithmetic function, a topic extensively studied by K.-H. \ Indlekofer, S.\ Wehmeier, and L.\ G.\ Lucht \cite{MR2116969}. In this section, we employ and modify their methods to establish the validity of \lemref{sumofvarthetap} and \lemref{sumofomegap}.

Consider an arithmetic function $F: \N \to \C$ and a strictly increasing sequence of natural numbers $\nseq \coloneqq (n_i)$. We define the mean behavior of $F$ along $\nseq$ as
\begin{align}\label{meanbehavior}
    \mathbb{M}(F,\nseq,x) \coloneqq \frac{1}{x} \sum_{i \leq x} F(n_i).
\end{align}
If $F$ is multiplicative, there exists a unique multiplicative function $G$ such that
\begin{align}\label{Dirichletinverse}
    F(n) = \sum_{d\mid n}G(d).
\end{align}
In fact, such $G$ is given by
$$G(n) \coloneqq \sum_{d\mid n} F(d) \mu(n/d),$$
where $\mu$ is the  M\"obius function. Substituting \eqref{Dirichletinverse} into \eqref{meanbehavior}, we obtain
\begin{align}\label{meanbehaviorwithDirichletinverse}
    \mathbb{M}(F,\nseq,x) = \frac{1}{x} \sum_{i \leq x} \sum_{d \mid n_i} G(d) = \frac{1}{x} \sum_{d=1}^\infty G(d) \sum_{\substack{i \leq x \\ d \mid n_i}}1.
\end{align}
Let $\mathfrak{p}'_{n,k}$ be the increasing sequence of shifted primes $p-1$ where $p \equiv k \pmod n$. Let us denote the $i$-th term of $\mathfrak{p}'_{n,k}$ by $p_i-1$. By \eqref{meanbehaviorwithDirichletinverse}, we have
\begin{align}\label{meanbehavior1}
    \pi(x;n,k) \cdot \mathbb{M}(F,\mathfrak{p}'_{n,k},\pi(x;n,k)) = \sum_{d=1}^\infty G(d) \sum_{\substack{i\leq \pi(x;n,k) \\ d \mid p_i-1}} 1 = \sum_{d=1}^\infty G(d) \sum_{\substack{p \leq x \\ p \equiv k (n) \\ p \equiv 1 (d)}}1.
\end{align}
Consider the function $f(n) \coloneqq \prod_{q\mid n} \left(1-\frac{1}{q(q^2-1)}\right)$, where $q$ runs over the prime factors of $n$. Note that $f$ is multiplicative. Also, from \eqref{definitionofvartheta}, we see that $\vartheta_p = f(p-1)$ for all primes $p$. Thus, the series in \lemref{sumofvarthetap} can be expressed as
\begin{align}\label{meanbehavior2}
    \sum_{\substack{p \leq x \\ p \equiv k(n)}} \vartheta_p = \sum_{\substack{p \leq x \\ p \equiv k(n)}} f(p-1) = \sum_{i \leq \pi(x;n,k)}f(p_i-1) = \pi(x;n,k) \cdot \mathbb{M}(f,\mathfrak{p}'_{n,k},\pi(x;n,k)).
\end{align}
From \eqref{definitionofg}, we find that
$$g(n) = \sum_{d \mid n} f(d) \mu(n/d),$$
so that $f(n) = \sum_{d\mid n}g(d)$. Combining \eqref{meanbehavior1} and \eqref{meanbehavior2}, we obtain
    $$\sum_{\substack{p \leq x \\ p \equiv k(n)}} \vartheta_p = \sum_{d=1}^\infty g(d) \sum_{\substack{p \leq x \\ p\equiv k(n) \\ p\equiv 1 (d) }} 1 = \sum_{d \text{ sq.free}} g(d) \sum_{\substack{p \leq x \\ p \equiv k (n) \\ p \equiv 1(d)}} 1 =: T_{n,k}(x),$$
since $g(d) = 0$ if and only if $d$ is not a square-free integer. Our goal is to determine the asymptotic formula for $T_{n,k}(x)$.

Let $S^*$ be the set of all square-free positive integers. Define
$$S_n^* \coloneqq \{s \in S^* : (n,s) = 1\}, \hs \hs \hs rS_n^* \coloneqq \{rs : s \in S_n^*\}.$$
Let $\text{rad}(n)$ be the radical of $n$, i.e., the product of all prime factors of $n$. Note that
$$S^* = \bigsqcup_{r \mid \text{rad}(n)} rS_n^*.$$
Define
\begin{equation}\label{definitionTnkrx}
    T_{n,k,r}(x) \coloneqq \sum_{d \in rS_n^*} g(d) \sum_{\substack{p \leq x \\ p \equiv k(n) \\ p \equiv 1(d)}}1.
\end{equation}
Then, we have 
\begin{equation}\label{Tnkx}
    T_{n,k}(x) = \sum_{r \mid \text{rad}(n)} T_{n,k,r}(x).
\end{equation}
Fix $r \mid \text{rad}(n)$. Suppose there exists a prime factor $q \mid r$ such that $k \not \equiv 1 \pmod q$. Then, the conditions $p\equiv k \pmod n$ and $p \equiv 1 \pmod d$ for $d \in rS_n^*$ contradict each other. In this case, we have $T_{n,k,r}(x) = 0$ for all $x$.

On the other hand, suppose that for all prime factors $q \mid r$, we have $k \equiv 1 \pmod q$. Consider $d = rd' \in rS_n^*$, so that $(d',n) = 1$. By the Chinese remainder theorem, we have
\begin{equation}\label{ChineseRemainderTheorem}
    \sum_{\substack{p \leq x \\ p \equiv k(n)\\ p \equiv 1 (d)}}1 = \sum_{\substack{p \leq x \\ p \equiv k(n)\\ p \equiv 1 (d')}}1 = \pi(x; nd', a(n,d')),
\end{equation}
where $a(n,d')$ is the unique integer between $1$ and $nd'$ such that
$$a(n,d') \equiv \begin{cases} k \pmod n \\ 1 \pmod{d'} \end{cases}.$$
From \eqref{definitionTnkrx} and \eqref{ChineseRemainderTheorem}, we obtain
\begin{equation}\label{Tnkrx}
    T_{n,k,r}(x) = \sum_{rd' \in rS_n^*} g(rd')  \sum_{\substack{p \leq x \\ p \equiv k(n) \\ p \equiv 1 (d')}}1 = g(r)\sum_{d' \in S_n^*}g(d') \pi(x;nd',a(n,d')).
\end{equation}
Hence, by \eqref{Tnkx} and \eqref{Tnkrx},
\begin{equation}\label{Tnkxreformulated1}
    T_{n,k}(x) = \sum_{\substack{r \mid \text{rad}(n) \\ \forall q \mid r, k \equiv 1 (q) }} g(r) \sum_{d \in S_n^*}g(d)\pi(x;nd,a(n,d)).
\end{equation}
One can check that given a square-free $r$, $k \equiv 1 \pmod r$ if and only if $k \equiv 1 \pmod q$ for all $q \mid r$. Since $g$ is multiplicative, we have
\begin{equation}\label{sumofg}
    \sum_{\substack{r \mid \rad(n) \\ \forall q \mid r , k \equiv 1 (q)}} g(r) = \sum_{\substack{r \mid \rad(n) \\ k \equiv 1 (r)}} g(r) = \prod_{\ell' \mid (n,k-1)} (1+g(\ell')).
\end{equation}
Therefore, by \eqref{Tnkxreformulated1} and \eqref{sumofg},
\begin{equation}\label{Tnkxreformulated}
    T_{n,k}(x) = \prod_{\ell' \mid (n,k-1)}(1+g(\ell')) \sum_{d\in S_n^*} g(d) \pi(x;nd,a(n,d))
\end{equation}
Hence, to prove \lemref{sumofvarthetap}, it suffices to determine the asymptotic formula for
\begin{equation}\label{Tnk1x}
    T_{n,k,1}(x) \coloneqq \sum_{d\in S_n^*} g(d) \pi(x;nd,a(n,d)).
\end{equation}

\subsection{Proof of \lemref{sumofvarthetap}}\label{s4.2}
We use the Bombieri-Vinogradov theorem to compute the asymptotic formula for $T_{n,k,1}(x)$.
\begin{Theorem}[Bombieri-Vinogradov theorem]\label{BombieriLemma}
Let 
$$E(y;d,l) \coloneqq \pi(y;d,l) - \frac{\li(y)}{\phi(d)} \hs \hs \text{ and } \hs \hs E(x,d) \coloneqq \max_{y \leq x}\max_{(l,d) = 1} |E(y;d,l)|.$$
Then, given a constant $D > 0$, there exists a constant $C = C(D) > 0$ such that
$$\sum_{d \leq x^{1/2}\log^{-C}x} E(x,d) = O_D(x\log^{-D}x).$$
\end{Theorem}
\begin{proof}
    See \cite[Theorem 4]{https://doi.org/10.1112/S0025579300005313}. (Also, see \cite[Theorem 2]{MR0197414} and \cite{MR0194397}.)
\end{proof}
Fix $D = K > 0$ and choose $C = C(K)  > 0$ according to \thmref{BombieriLemma}. Set $z = x^{1/2}\log^{-C}x$ and $z' = z/n$. Clearly, $\log z' \ll_K \log x$, since $z' \ll_K x$. Note that $\pi(x; nd,a(n,d)) = 0$ for $d > x$. Thus, we may truncate the series in $\eqref{Tnk1x}$ as
\begin{equation}\label{truncatedTnk1x}
    T_{n,k,1}(x) = \sum_{\substack{ d \in S_n^* \\ d \leq z'}} g(d)\pi(x;nd,a(n,d)) + \sum_{\substack{d \in S_n^* \\ z' < d \leq x}} g(d)\pi(x;nd,a(n,d)) =: S_1 + S_2.
\end{equation}
Using the trivial bound $\pi(x;nd,a(n,d)) \leq x/nd + 1$, we have
    \begin{align*}
        |S_2| &\leq \sum_{\substack{d \in S_n^* \\ z' < d \leq x}}|g(d)| \left(\frac{x}{nd}+1\right) \leq   \sum_{z'< d} |g(d)|\left(\frac{x}{nz'} + 1\right).
     \end{align*}
From \eqref{definitionofg} and the multiplicativity of $g$, we have $|g(d)| \leq 1/d^2$ holds for any integers $d$. Thus,
\begin{equation}\label{sumofgdisbounded}
   \sum_{z' < d} |g(d)| \leq \sum_{z' < d} \frac{1}{d^2} \ll \int_{z'}^\infty \frac{dt}{t^2} = \frac{1}{z'}.
\end{equation}
Hence,
\begin{equation}\label{S2bound}
    |S_2| \ll \frac{1}{z'}\left(\frac{x}{nz'} + 1\right) = n\log^{2C}x + nx^{-1/2}\log^C x \ll_K n\log^{2C}x.
\end{equation}
On the other hand, using the notation from \thmref{BombieriLemma}, we can express $S_1$ as
\begin{align}\label{S1}
    S_1 \coloneqq \sum_{\substack{d \in S_n^* \\ d \leq z'}} g(d) \pi(x;nd,a(n,d)) = \sum_{\substack{d \in S_n^* \\ d \leq z'}} g(d) \frac{\li(x)}{\phi(nd)} + \sum_{\substack{d \in S_n^* \\ d \leq z'}} g(d)E(x;nd,a(n,d)).
\end{align}
By the Cauchy-Schwarz inquality, \eqref{sumofgdisbounded}, and \thmref{BombieriLemma}, the last series can be bounded as
\begin{equation}\label{CauchySchwarz}
\begin{split}
      \left| \sum_{\substack{d \in S_n^* \\ d\leq z'}}g(d)E(x,nd, a(n,d))\right| &\leq \left( \sum_{d\leq z'}|g(d)|\right)^{1/2} \left( \sum_{d\leq z'}E(x,nd)\right)^{1/2} \\
    &\ll_K \left(nx^{-1/2}\log^C x\right)^{1/2} \left(x\log^{-K}x\right)^{1/2} = n^{1/2}x^{1/4}\log^{(C-K)/2}x.
\end{split}
\end{equation}
Now, focusing on the  main term of \eqref{S1}, since $\phi$ is multiplicative and $(n,d) = 1$, we have
\begin{align}\label{S1main}
    \sum_{\substack{d \in S_n^* \\ d \leq z'}} g(d) \frac{\li(x)}{\phi(nd)} = \frac{\li(x)}{\phi(n)}\sum_{\substack{d \in S_n^* \\ d \leq z'}} \frac{g(d)}{\phi(d)} = \frac{\li(x)}{\phi(n)}\sum_{d \in S_n^*} \frac{g(d)}{\phi(d)} + O\left( \frac{\li(x)}{\phi(n)}\sum_{z' < d} \frac{|g(d)|}{\phi(d)}\right)
\end{align}
Since $|g(d)|/\phi(d) \leq |g(d)|$, using the same argument from \eqref{sumofgdisbounded}, we can bound the error term in the above expression as follows:
\begin{align}\label{S1maintail}
  \frac{\li(x)}{\phi(n)} \sum_{z' < d} \frac{|g(d)|}{\phi(d)} \leq \frac{\li(x)}{\phi(n)} \sum_{z' < d} |g(d)|\ll \frac{\li(x)}{\phi(n)} n x^{-1/2}\log^Cx \ll nx^{1/2}\log^{C-1}x.
\end{align}
Since we are taking $n \leq \log^K x$, \eqref{S2bound}, \eqref{CauchySchwarz}, and \eqref{S1maintail} can be bounded by $O_{\delta,K}(x^{1/2+\delta})$ for any $\delta > 0$. Combining \eqref{truncatedTnk1x}, \eqref{S2bound}, \eqref{S1}, \eqref{S1main}, and \eqref{S1maintail}, we obtain
\begin{equation}\label{Tnk1xfinalized}
    T_{n,k,1}(x) = \frac{\li(x)}{\phi(n)}\sum_{d \in S_n^*}\frac{g(d)}{\phi(d)} + O_{\delta,K} \left(x^{1/2+\delta}\right),
\end{equation}
for sufficiently large $x$. Finally, note that $g/\phi$ is a multiplicative function. One can easily check that
\begin{align}\label{eulerproduct}
    \sum_{d\in S_n^*} \frac{g(d)}{\phi(d)} = \prod_{\ell \nmid n}\left(1 + \frac{g(\ell)}{\phi(\ell)}\right).
\end{align}
Thus, combining \eqref{Tnkxreformulated}, 
\eqref{Tnk1x}, \eqref{Tnk1xfinalized}, and \eqref{eulerproduct}, for any choice of $\delta > 0$, we obtain
$$T_{n,k}(x) = \frac{\li(x)}{\phi(n)}\prod_{\ell'\mid (n,k-1)} (1+g(\ell')) \prod_{\ell \nmid n} \left( 1 + \frac{g(\ell)}{\phi(\ell)}\right) + O_{\delta,K} \left(x^{1/2+\delta}\right),$$
for sufficiently large $x$. (The finite product of \eqref{Tnkxreformulated} is bounded by 1, and thus the error term remains unchanged.) This completes the proof of \lemref{sumofvarthetap}.
\subsection{Proof of \lemref{sumofomegap}}\label{s4.3}
We use the Siegel-Walfisz theorem to prove \lemref{sumofomegap}.
\begin{Theorem}[Siegel-Walfisz Theorem]\label{siegelwalfisz}
Let $D > 0$ and $q \leq \log^D x$. Then, there exists an ineffective constant $c$ that only depends on $D$ such that
$$\pi(x;q,a) = \frac{\li(x)}{\phi(q)} + O\left( x \exp \left( -c \sqrt{\log x}\right) \right).$$
\end{Theorem}
\begin{proof}
    See \cite[Hilfssatz 3]{MR1545584}
\end{proof}
Take $D = K + 1 > 1$. Fix positive integers $n$ and $m$ so that $\mathcal{M} \leq \log^{K+1} x$. Let $P(x;n,k)$ be the set of all primes $p \leq x$ for which $p \equiv k \pmod n$. Theis set can be partitioned as follows:
$$P(x;n,k) = \{p \in P(x;n,k) : (p,\underline{m}) = 1\} \cup \{p \in P(x;n,k) : (p,\underline{m}) \neq 1\} =: P_1(x) \cup P_2(x).$$
Note that $P_2(x)$ is a finite set. From the definition of $\underline{m}$, we observe that $(p,\underline{m}) \neq 1$ if and only if $p \mid m$. Thus, $\#P_2(x) \leq \log_2 m$, by the trivial bound. Since $\omega_p(m) \leq 1$ by definition, for any choice of $m$, we have
\begin{equation}\label{p2final}
    \sum_{p \in P_2(x)} \omega_p(m) \leq \log_2(m) =  O(\log m)
\end{equation}

Now, let us focus on the sum of $\omega_p(m)$ over $p \in P_1(x)$. As observed in \eqref{omegadepends}, the value of $\omega_p(m)$ is determined by the congruence class of $p$ modulo $\underline{m}$. Let us partition $P_1(x)$ as follows:
$$P_1(x) = \bigsqcup_{\substack{ 1 \leq a \leq \underline{m} \\ (a,\underline{m}) = 1}} \{p \in P_1(x) : p \equiv a  \mkern-12mu \pmod {\underline{m}}\}.$$
Thus, we have
\begin{equation}\label{last}
\sum_{p \in P_1(x)} \omega_p(m) = \sum_{\substack{1 \leq a \leq \underline{m} \\ (a,\underline{m}) = 1}} \sum_{\substack{p \leq x \\ p \equiv k (n) \\ p \equiv a (\underline{m})}} \omega_p(m) = \sum_{\substack{1 \leq a \leq \underline{m} \\ (a,\underline{m}) = 1}}\sum_{\substack{p \leq x \\ p \equiv k (n) \\ p \equiv a (\underline{m})}} \omega_a(m).\end{equation}

Consider a pair of integers $k$ and $a$, each of which is coprime to $n$ and $\underline{m}$, respectively. We now establish that
\begin{equation}\label{establish}
    \sum_{\substack{p \leq x \\ p \equiv k (n) \\ p \equiv a (\underline{m})}} \omega_a(m) = \sum_{\substack{1 \leq r \leq \mathcal{M}  \\ (r,\mathcal{M} ) = 1 \\ r \equiv k (n) \\ r \equiv a (\underline{m})}} \sum_{\substack{p \leq x \\ p \equiv r(\mathcal{M} )}}\omega_r(m) = \sum_{\substack{ 1 \leq r \leq \mathcal{M}  \\ (r,\mathcal{M} ) = 1 \\ r \equiv k (n) \\ r \equiv a (\underline{m})}}\omega_r(m) \pi(x;\mathcal{M} ,r).
\end{equation}
The latter equality is obvious. Let us validate the former one.

Suppose that the conditions $p \equiv k \pmod n$ and $p \equiv a \pmod {\underline{m}}$ are contradictory. In this case, no such $r$ exists in the interval $[1,\mathcal{M}]$ satisfying $r \equiv k \pmod n$ and $r \equiv a \pmod {\underline{m}}$. Hence, the sums are empty, and the equality holds. Suppose otherwise. Since $\mathcal{M} $ is chosen to be the least common multiple of $n$ and $\underline{m}$, there must be a unique integer $r$ in the interval $[1,\mathcal{M} ]$ for which $r \equiv k \pmod n$ and $r \equiv a \pmod {\underline{m}}$. By \remref{omegadepends}, we have $\omega_a(m) = \omega_r(m)$, and hence
$$\sum_{\substack{p \leq x \\ p \equiv k (n) \\ p \equiv a (\underline{m})}} \omega_a(m) = \sum_{\substack{p \leq x \\ p \equiv r (\mathcal{M} )}} \omega_a(m) = \sum_{\substack{p \leq x \\ p \equiv r (\mathcal{M} )}} \omega_r(m).$$
Hence, \eqref{establish} is verified.

From \eqref{last}, \eqref{establish}, and \thmref{siegelwalfisz}, we obtain
\begin{equation}\label{3.39}
    \begin{split}
        \sum_{p \in P_1(x)} \omega_p(m) = \sum_{\substack{1 \leq a \leq \underline{m} \\ (a,\underline{m}) = 1}} \sum_{\substack{1 \leq r \leq \mathcal{M}  \\ (r,\mathcal{M} ) = 1 \\ r \equiv k (n) \\ r \equiv a (\underline{m})}} \omega_r(m) \pi(x;\mathcal{M} ,r) &= \sum_{\substack{1 \leq r \leq \mathcal{M}  \\ (r,\mathcal{M} ) = 1 \\ r \equiv k (n)}} \omega_r(m) \pi(x;\mathcal{M} ,r) \\
        &=\sum_{\substack{1\leq r \leq \mathcal{M}  \\ (r,\mathcal{M} ) = 1 \\ r\equiv k(n)}} \omega_r(m)\left(\frac{\li(x)}{\phi(\mathcal{M} )} + O(x\exp\left(-c\sqrt{\log x}\right)\right).
    \end{split}
\end{equation}
Since $\#\{1 \leq r \leq \mathcal{M}  : (r,\mathcal{M} ) = 1, r \equiv k \pmod n\} \leq \phi(\underline{m}) \leq m$ and $\omega_r(m) \leq 1$, we have
\begin{equation}\label{3.40}
    \sum_{\substack{1 \leq r \leq \mathcal{M}  \\ (r,\mathcal{M} ) = 1 \\ r \equiv k (n)}}\omega_r(m) \leq m.
\end{equation}
From \eqref{p2final}, \eqref{3.39}, and \eqref{3.40}, we obtain
$$
 \sum_{\substack{p \leq x \\ p \equiv k(n)}} \omega_p(m) = O(\log m) +  \frac{\li(x)}{\phi(\mathcal{M} )}\sum_{\substack{1 \leq r \leq \mathcal{M}  \\ (r,\mathcal{M} ) =1 \\ r \equiv k(n)}} \omega_r(m)+ O\left( mx \exp \left( -c \sqrt{\log x}\right) \right).
$$
Since $O(\log m)$ can be absorbed into $O(mx \exp(-c\sqrt{\log x})$, \lemref{sumofomegap} is verified.
\section{Proofs of the Main Results}\label{s5}
\subsection{Proof of \propref{Densities}}\label{s5.1}
Let us prove the first claim. Fix $n > 2$. Recall the definition
$$C^{\Cyc}_{n,k} \coloneqq \frac{1}{\phi(n)} \prod_{\ell'\mid (n,k-1)}(1+g(\ell')) \prod_{\ell \nmid n}\left( 1 + \frac{g(\ell)}{\phi(\ell)}\right).$$
Let $n_k$ be the radical of $\gcd(n,k-1)$. For a fixed $n$, the density is determined by the number of prime factors of $n_k$. Since $-1 < g(\ell) < 0$ for any prime $\ell$, the more prime factors $n_k$ has, the smaller the density. Note that
\begin{align*}
    n_1 &= \rad(\gcd(n,0)) = \rad(n) \\
    n_{-1} &= \rad(\gcd(n,n-2)) = \begin{cases}
        1 & \text{ if } n \text{ is odd} \\
        2 & \text{ otherwise}
    \end{cases}.
\end{align*}
Thus, we have $n_{-1} \mid n_1$. Furthermore, for any $k$ coprime to $n$, it is clear from the definition that $n_{-1} \mid n_k \mid n_1$.

If $n$ is a power of two, one may note that $n_1 = n_{-1} = 2$. Hence, $C^{\Cyc}_{n,1} = C^{\Cyc}_{n,k} = C^{\Cyc}_{n,-1}$ for any odd $k$. If $n$ has an odd prime factor $q$, then $q \nmid n_{-1} = 1$, while $q \mid n_1 = \rad(n)$. Thus, we have $C^{\Cyc}_{n,1} < C^{\Cyc}_{n,-1}$.

Let us now check the second claim. We need to prove that $C^{\Div(m)}_{n,k}$ does not depend on the choice of $k$ if $(n,m) = 1$ or $m \in \{2,4\}$. First, assume that $n$ and $m$ are coprime. By the definition of $\underline{m}$, $n$ and $\underline{m}$ are also coprime, and hence $\mathcal{M} = n\underline{m}$. By the  Chinese remainder theorem, the natural projection
$$\{a \in (\Z/\mathcal{M}\Z)^\times : a \equiv k \mkern-12mu \pmod n\} \to (\Z/\underline{m}\Z)^\times,$$
is a bijection. Since $\omega_r(m)$ is determined by the congruence class of $r$ modulo $\underline{m}$, we obtain
$$\sum_{\substack{1 \leq r \leq \mathcal{M}\\ (r,\mathcal{M}) =1 \\ r \equiv k(n)}}\omega_r(m) = \sum_{\substack{1 \leq a \leq \underline{m} \\ (a,\underline{m}) = 1}} \omega_a(m) = \phi(\underline{m}) C^{\Div(m)}.$$
Therefore, we have
$$
C^{\Div(m)}_{n,k} \coloneqq \frac{1}{\phi(\mathcal{M})} \sum_{\substack{1 \leq r \leq \mathcal{M} \\ (r,\mathcal{M}) = 1 \\ r \equiv k (n)}} \omega_r(m) = \frac{\phi(\underline{m})C^{\Div(m)}}{\phi(\mathcal{M})} = \frac{C^{\Div(m)}}{\phi(n)},$$
and this quantity does not depend on the choice of $k$. 

Now, suppose $m$ is either 2 or 4. In both cases, we have $\underline{m} = 2$. It suffices to prove for the case that $n$ is even. In this case, we have $\mathcal{M}=n$, and hence
$$C^{\Div(m)}_{n,k} = \frac{1}{\phi(n)} \sum_{\substack{1 \leq r \leq n \\ (r,n) = 1 \\ r \equiv k (n)}} \omega_r(m) = \frac{\omega_k(m)}{\phi(n)}.$$
Since $k$ is coprime to $n$, $k$ must be odd. Thus, we have $\omega_k(2) = \omega_1(2)$ and $\omega_k(4) = \omega_1(4) = \omega_3(4)$ for any $k$ coprime to $n$. Each of these is equal to $C^{\Div(2)}$ and $C^{\Div(4)}$, respectively. Therefore, we have
$$C^{\Div(m)}_{n,k} = \frac{C^{\Div(m)}}{\phi(n)}$$
if $m \in \{2,4\}$. This completes the proof.
\subsection{Proofs of the Main Results}\label{s5.2}
In this section, we prove \thmref{maintheorem}.
\begin{Theorem}\label{maintheorem1}
Let $\epsilon > 0, x>0$, and $K > 0$. Let $A \coloneqq A(x)$ and $B \coloneqq B(x)$ be positive integers satisfying
$$x^{\epsilon} \leq A, B \leq x^{1-\epsilon}, \hs AB \geq x^{1+\epsilon}.$$
Fix $\sigma < \frac{\epsilon}{16+4\epsilon}$. Then, for any $n \leq \log^K x$ and $k$ coprime to $n$, we have
$$\frac{1}{4AB} \sum_{|a| \leq A } \sum_{|b| \leq B} \pi^{\Cyc}_{a,b}(x;n,k) = C^{\Cyc}_{n,k} \li(x) + O_{\epsilon,\sigma, K} \left(\frac{x^{1-\sigma}}{n}\right),$$
for sufficiently large $x$.
\end{Theorem}
\begin{proof}
From \lemref{NCnk0} and \lemref{sumofvarthetap}, we have
\begin{align*}
        N^{\Cyc}_{n,k}(A,B;x) &= 4AB \left( \frac{\li(x)}{\phi(n)} \prod_{\ell' \mid (n,k-1)}(1+g(\ell'))\prod_{\ell \nmid n}\left( 1 + \frac{g(\ell)}{\phi(\ell)}\right) + O_{\delta,K} \left( x^{1/2+\delta}\right) \right) + O_{\epsilon,\sigma,K}\left( \frac{1}{n}ABx^{1-\sigma}\right)\\
    &= 4AB \left( C^{\Cyc}_{n,k} \li(x) + O_{\delta,K} \left(x^{1/2 + \delta}\right) + O_{\epsilon,\sigma,K}\left(\frac{x^{1-\sigma}}{n}\right)\right),
\end{align*}
for sufficiently large $x$. Take $\delta = 1/3 - \sigma$, which is guaranteed to be positive, since $\sigma < 1/52$ by the given conditions. Then, we have
$$nx^{1/2 +\delta}\leq x^{5/6 - \sigma} \log ^Kx \ll_{K} x^{1-\sigma},$$
so $O_{\delta,K}(x^{1/2+\delta})$ can be absorbed into $O_{\epsilon,\sigma,K}(x^{1-\sigma}/n)$. Dividing both sides by $4AB$, we obtain the desired results.
\end{proof}

\begin{Theorem}\label{maintheorem2}
Let $\epsilon > 0$ and $x > 0$. Let $A \coloneqq A(x)$ and $B \coloneqq B(x)$ be positive integers satisfying
$$x^{\epsilon} \leq A, B \leq x^{1-\epsilon}, \hs AB \geq x^{1+\epsilon}.$$
Fix $K > 0$. Choose positive integers $n$ and $m$ such that $\mathcal{M} \leq \log^{K+1} x$, where $\mathcal{M}$ is defined in \eqref{definitionofM}. Then, we have
$$\frac{1}{4AB} \sum_{|a| \leq A } \sum_{|b| \leq B} \pi^{\Div(m)}_{a,b}(x;n,k) =  C^{\Div(m)}_{n,k}\li(x) + O_{\epsilon, K}\left(\frac{x}{n\log^{K+1} x} \right),$$
for sufficiently large $x$.
\end{Theorem}
\begin{proof}
We fix $0 < \epsilon < 1/3$. Otherwise, the conditions on $A$ and $B$ will be vacuous. We set $\tau = 1/4-11\epsilon/16 > 0$.

The condition on $\mathcal{M}$ implies that $n \leq \log^{K+1} x$ and $m \leq \log^{2K+2}x$. Therefore, for sufficiently large $x$ sufficiently large, we have $n \leq x^\tau$. By \lemref{NDnk0} and \lemref{sumofomegap},  there exists $c > 0$ that only depends on $K$ such that
\begin{align*}
    N^{\Div(m)}_{n,k}(A,B;x) &= 4AB\left(\frac{\li(x)}{\phi(\mathcal{M})} \sum_{\substack{1\leq r \leq \mathcal{M}  \\ (r,\mathcal{M} ) = 1 \\ r \equiv k (n)}} \omega_r(m) + O\left(mx\exp(-c\sqrt{\log x})\right) \right) + O_{\epsilon, K}\left( \frac{ABx}{n\log^{K+1}x}\right) \\
    &=4AB\left(C^{\Div(m)}_{n,k}\li(x) +O\left(mx\exp(-c\sqrt{\log x})\right) + O_{\epsilon,K}\left(\frac{x}{n\log^{K+1}x}\right)\right),
\end{align*}
for sufficiently large $x$. Since $n \leq \log^{K+1} x$ and $m \leq \log^{2K+2}x$, we have 
\begin{align*}
    &\hs nm \log^{K+1}x \cdot \exp\left(-c\sqrt{\log x}\right) \leq \log^{4K+4} x \cdot  \exp \left(-c\sqrt{\log x}\right) = o_K(1) \\
    \Longrightarrow &\hs mx\exp\left(-c\sqrt{\log x}\right) = o_K\left( \frac{x}{n\log^{K+1}x}\right).
\end{align*}
Thus, $O(mx\exp(-c\sqrt{\log x})$ can be absorbed into $O_{\epsilon, K} (x/n\log^{K+1}x)$. Dividing both sides by $4AB$, we obtain the desired results.
\end{proof}

\section{Tables of Values of \texorpdfstring{$C^{\Cyc}_{n,k}$}{CCnk} and \texorpdfstring{$C^{\Div(m)}_{n,k}$}{CDmnk}}\label{s6}
\subsection{Average cyclicity density for primes in an arithmetic progression}\label{s6.1}
Let $n$ be a power of a prime $p$ and let $k$ be coprime to $n$. Recall the definition of $n_k$, the radical of $\gcd(n,k-1)$. If $p = 2$, then for any odd $k$, we have $n_k = 2$, meaning there is a single possible value for $C^{\Cyc}_{n,k}$. If $p$ is an odd prime, then $n_k$ can be either 1 or $p$, depending on whether $k \equiv 1 \pmod p$. Thus, there are exactly two possible values of $C^{\Cyc}_{n,k}$. The values of $\phi(n)C^{\Cyc}_{n,k}$ for small prime powers are listed in Table \tabref{Table1}. (The author used $\phi(n)C^{\Cyc}_{n,k}$ instead of $C^{\Cyc}_{n,k}$ to facilitate comparison between the densities of different moduli.) Notice that the quantity in the cell for $p = 2$ and $k \equiv 1 \pmod p$ equals the average cyclicity density $C^{\Cyc}$ determined by Banks and Shparlinski. 

As seen in Table \tabref{Table1}, as $p$ increases, the difference between the two densities along each column gets smaller. This is because their ratio is $1+g(p) = 1-\frac{1}{p(p^2-1)}$, which approaches to $1$ as $p \to \infty$.

Now, suppose $n$ is supported on two distinct odd primes $p$ and $q$. In this case, there are four possible values of $n_k$: $1, p, q,$ and $pq$. Table \tabref{Table2} lists the values of $\phi(n)C^{\Cyc}_{n,k}$ for $n$ supported on two small odd primes. The table shows that the densities along each column are distinct. Observing this, one pose the following question.
\begin{Question}\label{q1}
    For a fixed $n$ supported on $s$ distinct odd primes, are there $2^s$ distinct values of $C^{\Cyc}_{n,k}$ as $k$ varies?
\end{Question}
Recall that
$$g(\ell) = -\frac{1}{\ell(\ell^2-1)}.$$
Define
$$h(N) \coloneqq \prod_{\ell \mid N} \left( 1 - \frac{1}{\ell(\ell^2-1)}\right),$$
where the product is taken over all prime factors of $N$. Then, \questref{q1} is equivalent to determining whether there exist distinct square-free integers $N$ and $M$ for which $h(N) = h(M)$. Moreover, since $h$ is multiplicative, we can assume that $N$ and $M$ are coprime.

Suppose no such integer pair of integers exist. Then, for any $n$ supported on $s$ odd primes, there are $2^s$ distinct values of $C^{\Cyc}_{n,k}$. However, if such a pair $N$ and $M$ exists, then one can select $n = NM$ and $k_1, k_2$ coprime to $n$ for which $n_{k_1} = N$ and $n_{k_2} = M$, yielding $C^{\Cyc}_{n,k_1} = C^{\Cyc}_{n,k_2}$.

This question is connected to the following group theory problem. Note that $\ell(\ell^2-1) = |\SL_2(\Z/\ell\Z)|$. Say $a$ and $b$ are positive square-free integers with $a \mid b$. Define
$$S(a|b) := \{ g \in \SL_2(\Z/b\Z) : g \not \equiv I \mkern-12mu \pmod p \text{ for any } p \mid a\},$$
where $I$ denotes the identity element. By the Chinese remainder theorem, we have
$$|S(a|b)| =|\SL_2(\Z/(b/a)\Z)|\prod_{\ell \mid a} \left( |\SL_2(\Z/\ell\Z)| - 1\right) = |\SL_2(\Z/(b/a)\Z)|\prod_{\ell \mid a} (\ell^3-\ell - 1).$$
Suppose that such coprime square-free integer pairs $N$ and $M$ exist with $h(N) = h(M)$. Then,
\begin{align*}
    h(N) = h(M) \hs &\Longleftrightarrow \hs |\SL_2(\Z/NM\Z)| h(N) = |\SL_2(\Z/NM\Z)| h(M) \\
    \hs &\Longleftrightarrow \hs |\SL_2(\Z/M\Z)| \prod_{\ell \mid N}(\ell^3-\ell-1) = |\SL_2(\Z/N\Z)| \prod_{\ell' \mid M} (\ell'^3-\ell'-1) \\
    \hs &\Longleftrightarrow \hs |S(M|NM)| = |S(N |NM)|.
\end{align*}
Therefore, \questref{q1} can be restated as follows:
\begin{Question}
    Does there exists a positive square-free integer $n$ such that
    $$\#\{g \in \SL_2(\Z/n\Z) : g \not \equiv I \mkern-12mu \pmod p \text{ for all } p \mid m_1\} = \#\{g \in \SL_2(\Z/n\Z) : g \not \equiv I \mkern-12mu \pmod q \text{ for all } q \mid m_2\},$$
    for some distinct coprime factors $m_1$ and $m_2$ of $n$?
\end{Question}
The author confirmed using Magma, a computer-assisted calculating software, that no such integer pairs $N$ and $M$ exist in the range $[1,10000]$.

\subsection{Average divisibility density for primes in arithmetic progression}\label{s6.2}
The values of $C^{\Div(m)}_{n,k}$ are listed in Table \tabref{Table3} for $n, m \in \{2,3,4,5,6,8\}$. (Note that $C^{\Div(m)}_{2,1}$ is equal to $C^{\Div(m)}$ by definition.)

\propref{Densities} states that there is no average congruence class bias when $(m,n) = 1$ or $m \in \{2,4\}$. However, these are not the only necessary conditions; there are cases where $(n,m) \neq 1$ and $m \not \in \{2,4\}$, yet $C^{\Div(m)}_{n,k} = C^{\Div(m)}/\phi(n)$. For the readers' convenience, such values are marked with $\dagger$ signs in Table \tabref{Table3}.

As mentioned in \sectionref{s1.1}, Cojocaru computed that the individual density $C^{\Div(m)}_E$ is approximately $1/m$. We now prove that the average density $C^{\Div(m)}$ is bounded below by $1/m$, consistent to the individual case.
\begin{Lemma}\label{qlemma}
 Consider a positive integer $m$. Let $1 \leq r \leq m$ with $r$ coprime to $m$. Then, we have
 $$\omega_r(m) \geq \frac{1}{m}.$$
\end{Lemma}
\begin{proof}
    Since $\omega_r$ is a multiplicative function, it suffices to prove for the case of $m = q^j$. (The definition of $\omega_r$ can be found from \eqref{omegardefine}.) Say $r \not \equiv 1 \pmod {q^{\lceil j/2\rceil}}$. Then, we have
    $$\omega_r(q^j) = \frac{1}{q^{j-1}(q-1)} \geq \frac{1}{q^j}.$$
    Now, suppose $r \equiv 1 \pmod {q^{\lceil j/2\rceil}}$. Suppose $j = 2j' + 1$ for some non-negative integer $j'$. Thus, we have
    \begin{align*}
        q(q^{j'+1}+q^{j'}-1) \geq q^{j'}(q^2-1) \hs \implies \hs \frac{q(q^{j'+1} + q^{j'} - 1)}{q^{j'}(q^2-1)}\cdot \frac{1}{q^{2j'+1}} \geq \frac{1}{q^{2j'+1}} \hs \implies \hs \omega_r(q^j) \geq \frac{1}{q^j}.
    \end{align*}
    On the other hand, suppose $j = 2j'$ for some positive integer $j'$. Then, we have
    $$q^{j'+1} + q^{j'} -1 \geq q^{j'-1}(q^2-1) \hs \Longrightarrow \hs \frac{q^{j'+1} + q^{j'} - 1}{q^{j'-1} (q^2-1)} \cdot \frac{1}{q^{2j'}} \geq \frac{1}{q^{2j'}} \hs \Longrightarrow \hs \omega_r(q^j) \geq \frac{1}{q^j}.$$
    This completes the proof.
\end{proof}

\begin{Proposition}
    For any positive $m$, we have
    $$C^{\Div(m)} \geq \frac{1}{m}.$$
\end{Proposition}
\begin{proof}
    Recall the definition of $C^{\Div(m)}:$
    $$C^{\Div(m)} \coloneqq \frac{1}{\phi(\underline{m})} \sum_{\substack{ 1 \leq k \leq \underline{m} \\ (k,\underline{m}) = 1}} \omega_k(m) \geq \frac{1}{\phi(\underline{m})}\sum_{\substack{ 1 \leq k \leq \underline{m} \\ (k,\underline{m}) = 1}} \frac{1}{m} = \frac{1}{m},$$
    by \lemref{qlemma}. This completes the proof.
\end{proof}

\clearpage

\begin{table}[H]
\centering
\begin{tabular}{|c|c|c|c|c|c|} \hline
& & & & & \\
$p$ & 3 & 5 & 7 & 11 & 13 \\
& & & & & \\ \hline 
& & & & & \\
$k \equiv 1 \pmod p$ & $0.79644\ldots$ & $0.80866\ldots$ & $0.81173\ldots$ & $0.81320\ldots$ & $0.81341\ldots$ \\
& & & & & \\ \hline
& & & & & \\
$k \not \equiv 1 \pmod p$ & $0.83107\ldots$ & $0.81545\ldots$ & $0.81415\ldots$ & $0.81381\ldots$ & $0.81378\ldots$ \\
& & & & & \\ \hline
\end{tabular}
\caption{Values of $\phi(n)C^{\Cyc}_{n,k}$ for $n$ an odd prime power} \label{Table1}
\end{table}

\begin{table}[H]
\begin{center}
\begin{tabular}{|c|c|c|c|c|c|c|} \hline
& & & & & & \\
$(p,q)$ & (3,5) & (3,7) & (3,11) & (5,7) & (5,11) & (7,11) \\
& & & & & & \\ \hline 
& & & & & & \\
$\begin{matrix}
    k \equiv 1 \pmod p \\
    k \equiv 1 \pmod q
\end{matrix}$ & $0.79145\ldots$ & $0.79446 \ldots$ & $0.79589\ldots$ & $0.80665\ldots$ & $0.80810 \ldots$ & $0.81118\ldots$ \\
& & & & & & \\ \hline
& & & & & & \\
$\begin{matrix}
    k \equiv 1 \pmod p \\
    k \not \equiv 1 \pmod q
\end{matrix}$ & $0.79810\ldots$ & $0.79683\ldots$ & $0.79650 \ldots$ & $0.80906\ldots$ & $0.80872\ldots$ & $0.81179\ldots$ \\
& & & & & & \\ \hline
& & & & & & \\
$\begin{matrix}
    k \not \equiv 1 \pmod p \\
    k \equiv 1 \pmod q
\end{matrix}$ & $0.82586\ldots$ & $0.82900\ldots$ & $0.83050\ldots$ & $0.81343\ldots$ & $0.81490\ldots$ & $0.81360\ldots$ \\
& & & & & & \\ \hline
& & & & & & \\
$\begin{matrix}
    k \not \equiv 1 \pmod p \\
    k \not \equiv 1 \pmod q
\end{matrix}$ & $0.83280\ldots$ & $0.83148 \ldots$ & $0.83113\ldots$ & $0.81586\ldots$ & $0.81551\ldots$ & $0.81422\ldots$ \\
& & & & & & \\ \hline
\end{tabular}
\end{center}
\caption{Values of $\phi(n)C^{\Cyc}_{n,k}$ for $n$ a product of two prime powers}
\label{Table2}
\end{table}
\clearpage

\begin{table}[H]
    \begin{center}
        \begin{tabular}{|c|c|c|c|c|c|c|} \hline \xrowht{20pt}
     \: $\displaystyle (n,k)$ \: & \: $m=2$ \: &\: $m=3$ \: &\: $m=4$ \:&\: $m=5$ \:& \: $m=6$ \: &\: $m=8$ \: \\ \hline \hline \xrowht{20pt}
     $(2,1) $ & $2/3$ & $7/16$ & $5/12$ & $23/96$ & $7/24$ & $11/48$ \\ \hline \hline \xrowht{20pt}
     $(3,1)$ & $1/3$ & $3/16$ & $5/24$ & $23/192$ & $1/8$ & $11/96$ \\ \hline \xrowht{20pt}
     $(3,2)$& $1/3$ & $1/4$ & $5/24$ & $23/192$ & $1/6$ & $11/96$ \\ \hline \hline \xrowht{20pt}
     $(4,1)$ & $1/3$ & $7/32$ & $5/24$ & $23/192$ & $7/48^\dagger$ & $5/48$ \\ \hline \xrowht{20pt}
     $(4,3)$ & $1/3$ & $7/32$ & $5/24$ & $23/192$ & $7/48^\dagger$ & $1/8$ \\ \hline \hline \xrowht{20pt}
     $(5,1)$ & $1/6$ & $7/64$ & $5/48$ & $5/96$ & $7/96$ & $11/192$ \\ \hline \xrowht{20pt}
     $(5,2)$ & $1/6$ & $7/64$ & $5/48$ & $1/16$ & $7/96$ & $11/192$ \\ \hline \xrowht{20pt}
     $(5,3)$ & $1/6$ & $7/64$ & $5/48$ & $1/16$ & $7/96$ & $11/192$ \\ \hline \xrowht{20pt}
     $(5,4)$ & $1/6$ & $7/64$ & $5/48$ & $1/16$ & $7/96$ & $11/192$ \\ \hline \hline \xrowht{20pt}
     $(6,1)$ & $1/3$ & $3/16$ & $5/24$ & $23/192$ & $1/8$ & $11/96^\dagger$ \\ \hline \xrowht{20pt}
     $(6,5)$ & $1/3$ & $1/4$ & $5/24$ & $23/192$ & $1/6$ & $11/96^\dagger$ \\ \hline \hline \xrowht{20pt}
     $(8,1)$ & $1/6$ & $7/64$ & $5/48$ & $23/384$ & $7/96^\dagger$ & $5/96$ \\ \hline \xrowht{20pt}
     $(8,3)$ & $1/6$ & $7/64$ & $5/48$ & $23/384$ & $7/96^\dagger$ & $1/16$ \\ \hline \xrowht{20pt}
     $(8,5)$ & $1/6$ & $7/64$ & $5/48$ & $23/384$ & $7/96^\dagger$ & $5/96$ \\ \hline \xrowht{20pt}
     $(8,7)$ & $1/6$ & $7/64$ & $5/48$ & $23/384$ & $7/96^\dagger$ & $1/16$ \\ \hline
    \end{tabular} 
    \end{center}
    \caption{Values of $C^{\Div(m)}_{n,k}$ for $n$ and $m$ small integers}
    \label{Table3}
\end{table}

\clearpage

\bibliographystyle{amsplain}
\bibliography{NewReferences}

\providecommand{\bysame}{\leavevmode\hbox to3em{\hrulefill}\thinspace}
\providecommand{\MR}{\relax\ifhmode\unskip\space\fi MR }
\providecommand{\MRhref}[2]{%
  \href{http://www.ams.org/mathscinet-getitem?mr=#1}{#2}
}
\providecommand{\href}[2]{#2}
\begin{thebibliography}{10}

\bibitem{MR4504664}
Yildirim Akbal and Ahmet~M. G\"{u}lo\u{g}lu, \emph{Cyclicity of elliptic curves modulo primes in arithmetic progressions}, Canad. J. Math. \textbf{74} (2022), no.~5, 1277--1309. \MR{4504664}

\bibitem{MR2376806}
Stephan Baier, \emph{The {L}ang-{T}rotter conjecture on average}, J. Ramanujan Math. Soc. \textbf{22} (2007), no.~4, 299--314. \MR{2376806}

\bibitem{MR2570668}
William~D. Banks and Igor~E. Shparlinski, \emph{Sato-{T}ate, cyclicity, and divisibility statistics on average for elliptic curves of small height}, Israel J. Math. \textbf{173} (2009), 253--277. \MR{2570668}

\bibitem{https://doi.org/10.1112/S0025579300005313}
E.~Bombieri, \emph{On the large sieve}, Mathematika \textbf{12} (1965), no.~2, 201--225.

\bibitem{Brau}
Julio Brau, \emph{Galois representations of elliptic curves and abelian entanglements}, 2015, Doctoral Thesis, Available at \url{https://scholarlypublications.universiteitleiden.nl/handle/1887/37019}.

\bibitem{brau2017character}
\bysame, \emph{Character sums for elliptic curve densities}, 2017, arXiv:1703.04154, Available at \url{https://arxiv.org/abs/1703.04154}.

\bibitem{campagna2022cyclic}
Francesco Campagna and Peter Stevenhagen, \emph{Cyclic reduction densities for elliptic curves}, Res. Number Theory \textbf{9} (2023), no.~3, Paper No. 61, 21. \MR{4623047}

\bibitem{MR1932460}
Alina~Carmen Cojocaru, \emph{On the cyclicity of the group of {$\mathbb{F}_p$}-rational points of non-{CM} elliptic curves}, J. Number Theory \textbf{96} (2002), no.~2, 335--350. \MR{1932460}

\bibitem{MR2076566}
\bysame, \emph{Questions about the reductions modulo primes of an elliptic curve}, Number theory, CRM Proc. Lecture Notes, vol.~36, Amer. Math. Soc., Providence, RI, 2004, pp.~61--79. \MR{2076566}

\bibitem{MR4033729}
\bysame, \emph{Primes, elliptic curves and cyclic groups}, Analytic methods in arithmetic geometry, Contemp. Math., vol. 740, Amer. Math. Soc., [Providence], RI, [2019] \copyright 2019, With an appendix by Cojocaru, Matthew Fitzpatrick, Thomas Insley and Hakan Yilmaz, pp.~1--69. \MR{4033729}

\bibitem{MR1677267}
Chantal David and Francesco Pappalardi, \emph{Average {F}robenius distributions of elliptic curves}, Internat. Math. Res. Notices (1999), no.~4, 165--183. \MR{1677267}

\bibitem{MR2235360}
M.~Z. Garaev, \emph{Character sums in short intervals and the multiplication table modulo a large prime}, Monatsh. Math. \textbf{148} (2006), no.~2, 127--138. \MR{2235360}

\bibitem{MR1055716}
Rajiv Gupta and M.~Ram Murty, \emph{Cyclicity and generation of points mod {$p$} on elliptic curves}, Invent. Math. \textbf{101} (1990), no.~1, 225--235. \MR{1055716}

\bibitem{MR207630}
Christopher Hooley, \emph{On {A}rtin's conjecture}, J. Reine Angew. Math. \textbf{225} (1967), 209--220. \MR{207630}

\bibitem{MR1204781}
Everett~W. Howe, \emph{On the group orders of elliptic curves over finite fields}, Compositio Math. \textbf{85} (1993), no.~2, 229--247. \MR{1204781}

\bibitem{MR2116969}
K.-H. Indlekofer, S.~Wehmeier, and L.~G. Lucht, \emph{Mean behaviour and distribution properties of multiplicative functions}, Comput. Math. Appl. \textbf{48} (2004), no.~12, 1947--1971. \MR{2116969}

\bibitem{JKW}
Seraphim Jarov, Alex Khadra, and Nahid Walji, \emph{Congruence class bias and the {L}ang-{T}rotter conjecture for families of elliptic curves}, Involve \textbf{17} (2024), no.~4, 569--592. \MR{4810094}

\bibitem{MR2563740}
Nathan Jones, \emph{Almost all elliptic curves are {S}erre curves}, Trans. Amer. Math. Soc. \textbf{362} (2010), no.~3, 1547--1570. \MR{2563740}

\bibitem{jones2023acyclicity}
Nathan Jones and Sung~Min Lee, \emph{On the acyclicity of reductions of elliptic curves modulo primes in arithmetic progressions}, 2023, arXiv:2206.00872, Available at \url{https://arxiv.org/abs/2206.00872}.

\bibitem{jones2021elliptic}
Nathan Jones and Kevin Vissuet, \emph{Elliptic curves with missing frobenius traces}, 2021, arXiv:2108.08727, Available at \url{https://arxiv.org/abs/2108.08727}.

\bibitem{lee2024opposingaveragecongruenceclass}
Sung~Min Lee, Jacob Mayle, and Tian Wang, \emph{Opposing average congruence class biases in the cyclicity and koblitz conjectures for elliptic curves}, 2024, arXiv:2408.16641, Available at \url{https://arxiv.org/abs/2408.16641}.

\bibitem{lmfdb}
The {LMFDB Collaboration}, \emph{The {L}-functions and modular forms database}, \url{https://www.lmfdb.org}, 2024, [Online; accessed 6 March 2024].

\bibitem{mayle2023serre}
Jacob Mayle and Rakvi, \emph{Serre curves relative to obstructions modulo 2}, Lu{C}a{NT}: {LMFDB}, computation, and number theory, Contemp. Math., vol. 796, Amer. Math. Soc., [Providence], RI, [2024] \copyright 2024, pp.~103--128. \MR{4732685}

\bibitem{MR0698163}
M.~Ram Murty, \emph{On {A}rtin's conjecture}, J. Number Theory \textbf{16} (1983), no.~2, 147--168. \MR{698163}

\bibitem{MR2470688}
Richard Taylor, \emph{Automorphy for some {$l$}-adic lifts of automorphic mod {$l$} {G}alois representations. {II}}, Publ. Math. Inst. Hautes \'{E}tudes Sci. (2008), no.~108, 183--239. \MR{2470688}

\bibitem{MR0197414}
A.~I. Vinogradov, \emph{The density hypothesis for {D}irichet {$L$}-series}, Izv. Akad. Nauk SSSR Ser. Mat. \textbf{29} (1965), 903--934. \MR{197414}

\bibitem{MR0194397}
\bysame, \emph{Correction to the paper of {A}. {I}. {V}inogradov ``{O}n the density hypothesis for the {D}irichlet {$L$}-series''}, Izv. Akad. Nauk SSSR Ser. Mat. \textbf{30} (1966), 719--720. \MR{194397}

\bibitem{MR1667099}
S.~G. Vl\v{a}du\c{t}, \emph{Cyclicity statistics for elliptic curves over finite fields}, Finite Fields Appl. \textbf{5} (1999), no.~1, 13--25. \MR{1667099}

\bibitem{MR1545584}
Arnold Walfisz, \emph{Zur additiven {Z}ahlentheorie. {II}.}, Math. Z. \textbf{40} (1936), no.~1, 592--607. \MR{1545584}

\bibitem{MR4765791}
Peng-Jie Wong, \emph{Cyclicity and exponent of elliptic curves modulo {$p$} in arithmetic progressions}, Q. J. Math. \textbf{75} (2024), no.~2, 757--777. \MR{4765791}

\bibitem{zywina2015possible}
David Zywina, \emph{On the possible images of the mod $\ell$ representations associated to elliptic curves over $\mathbb{Q}$}, 2015, arXiv:1508.07660, Available at \url{https://arxiv.org/abs/1508.07660}.

\end{thebibliography}

\end{document}